\documentclass[11pt]{article}
\usepackage{preamble}

\makeatletter
\renewcommand{\@fnsymbol}[1]{\ifcase#1\or 1\or 2\or 3\or 4\or 5\or 6\fi}
\makeatother

\usepackage{natbib}
 \bibpunct[, ]{(}{)}{,}{a}{}{,}%

\author{Sumin Kang\thanks{\href{mailto:suminkang@vt.edu}{suminkang@vt.edu}}~} 
\author{Manish Bansal\thanks{\href{mailto:bansal@vt.edu}{bansal@vt.edu}}}
\affil{Grado Department of Industrial and Systems Engineering, Virginia Tech}
\title{
    Bi-Parameterized Two-Stage Stochastic Min-Max and Min-Min Mixed Integer Programs
}
\date{\vspace{-3em}}

\begin{document}

\maketitle

\begin{abstract}
We introduce two-stage stochastic min-max and min-min integer programs with bi-parameterized recourse (BTSPs), where the first-stage decisions affect both the objective function and the feasible region of the second-stage problem. To solve these programs efficiently, we introduce Lagrangian-integrated $L$-shaped ($L^2$) methods, which guarantee exact solutions when the first-stage decisions are pure binary. For mixed-binary first-stage programs, we present a regularization-augmented variant of this method. 
Our computational results for a stochastic network interdiction problem show that the $L^2$ method outperforms a benchmark method, solving all instances in 23 seconds on average, while the benchmark method failed to solve any instance within 3600 seconds. The $L^2$ method also achieves optimal solutions, on average, 18.4 times faster for a stochastic facility location problem. Furthermore, we show that the $L^2$ method can effectively address distributionally robust optimization problems with decision-dependent ambiguity sets that may be empty for some first-stage decisions, achieving optimal solutions, on average, 5.3 times faster than existing methods.
\end{abstract}

\paragraph{Key words:} Stochastic integer programs, Bi-parameterized recourse, Interdiction problems, Decision-dependent uncertainty, Distributionally robust optimization

\section{Introduction} \label{sec:introduction}

Two-stage stochastic programming is a well-known framework for modeling decision-making under uncertainty, where decisions are made sequentially over two stages: an initial set of (first-stage) decisions before the uncertainty is revealed, followed by a set of recourse (second-stage) decisions that adapt to revealed outcomes. This framework has been used for a wide variety of applications such as network interdiction \citep{SongSmithSurvey}, healthcare \citep{yoon2021}, power systems \citep{zheng2013UnitCommit}, airline crew scheduling \citep{yenBirge2006}, wildfire planning \citep{ntaimo2012}, etc. In this framework, it is typically assumed that the first-stage decisions affect the right-hand side of the constraints in the second-stage problem where the recourse decisions are made. In this paper, we investigate \textit{bi-parameterized two-stage stochastic programming}, an extension of the conventional two-stage stochastic programming framework where the first-stage decisions affect the constraints and also the objective function of the second-stage problem. The formulation of bi-parameterized two-stage stochastic programs (BTSPs) is given by 
\begin{equation} \label{eq:btsp}
    \min_{x \in X} \Big\{f(x) + \E\big[Q(x, \om)\big]\Big\},
\end{equation}
where vector $x$ denotes the set of first-stage decision variables, Let $(q, r, G, W, T)$ be random data parameters associated with the second stage. For each scenario $\om \in \Om$, their realization, denoted by $(q(\om), r(\om), G(\om), W(\om), T(\om))$, occurs with probability $p(\om)$.
The recourse function $Q(x, \om)$ is defined for each $\om \in \Om$ as follows:
\begin{equation} \label{eq:Q}
\begin{aligned}
    Q(x, \om) := & \minmax \Big\{ q(\om) ^\top y + x^\top G(\om) y : W(\om) y = r(\om) - T(\om) x, \ y \in \Y \Big\}.
\end{aligned}
\end{equation}
The notation ``$\minmax$'' in \eqref{eq:Q} indicates that the recourse problem can be either a minimization problem or a maximization problem. Throughout the paper, we refer to BTSP~\cref{eq:btsp} as \textit{min-min problem} when the recourse problem~\cref{eq:Q} is a minimization, and \textit{min-max problem} when the recourse problem is a maximization problem. It can be easily seen that this formulation reduces to the conventional ``single-parameterized" two-stage stochastic program when $\minmax$ replaced with $\min$ and $G(\om) = 0$ for all $\om \in \Om$.

The first-stage feasible set is defined as $X := \{x \in \X: Ax = b\}$, where $\X \subseteq \R^{n_1}_+$ represents integrality restrictions on $x$, $A \in \Qu^{m_1 \times n_1}$, and $ b\in\Qu^{m_1}$. The function $f:X \to \R$ is a linear function that represents the first-stage objective. 
Let $Y(x, \om) := \{y \in \Y: W(\om) y = r(\om) - T(\om) x\}$ for $\om \in \Om$ where $\Y \subseteq \R^{n_2}_+$ represents integrality restrictions on $y$. Here, $q(\om)\in\Qu^{n_2}$, $G(\om) \in \Qu^{n_1 \times n_2}$, $T(\om) \in \Qu^{m_2 \times n_1}$, $W(\om) \in \Qu^{m_2 \times n_2}$, and $r(\om) \in \Qu^{m_2}$. 

A single-parameterized reformulation of BTSPs can be derived by introducing a proxy variable with constraints $\theta = x^\top G(\om) y$. However, since the resulting constraints depend on the parameterized coefficient, i.e., $x^\top G(\om)$ for $y$, the problem structure differs from the classical formulation of two-stage stochastic programming where coefficients associated with $y$ are independent of $x$.

The bi-parameterization introduces significant computational challenges, mainly due to the nonconvexity of the recourse function, even without integrality restrictions on the second-stage variables $y$. To illustrate this, consider the following example.

\begin{figure}[t]
    \centering
    \includegraphics[width=0.5\textwidth]{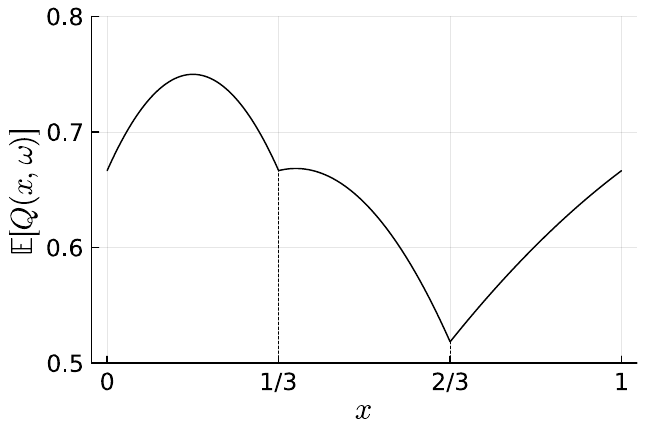}
    \caption{An example of $\E[Q(x, \om)]$.} 
    \label{fig:ex-nonconvex-recourse}
\end{figure}

\begin{example}
Let $\Om = \{\frac{1}{3}, \frac{2}{3}, 1\}$ and $p(\om) = 1/3$ for $\om \in \Om$. Consider the following recourse function:
% \begin{align}
$
    Q(x, \om) = \min \{(1+x)y_1 + (1+3x)y_2 : 
    y_1 \geq x - \om, \ y_2 \geq \om - x, \ y \in \R^2_+ \},
$
% \end{align}
 for $\om \in \Om$. The expected recourse function, i.e., $\E[Q(x, \om)] = \frac{1}{3}Q(x, \frac{1}{3}) + \frac{1}{3}Q(x, \frac{2}{3}) + \frac{1}{3}Q(x, 1)$, on $[0, 1]$ is nonconvex and nonsmooth, as illustrated in \Cref{fig:ex-nonconvex-recourse}.
We also note that, when integrality restrictions are imposed on $y$, the recourse function becomes discontinuous. 
$\blacksquare$\end{example}

These properties of the recourse function limit the applicability of traditional decomposition methods. For example, consider Benders decomposition (or the $L$-shaped method), which relies on \textit{optimality cuts} of the form $Q(x, \om) \geq \bar{\pi}^\top x +\bar{\pi}^0$ for all $x \in X$ to approximate the recourse function. For a single-parameterized recourse, we can solve the linear programming (LP) dual of the continuous relaxation of the recourse problem to obtain dual solutions and construct optimality cuts; while these cuts may not be tight, they are valid. However, when applied to BTSPs, such cuts are not valid for the recourse function due to its nonconvexity, making classical Benders decomposition methods inapplicable.
Alternatively, the dual decomposition method~\citep{Caroe1999} can be applied to BTSPs. In this approach, copies of $x$, denoted by $x(\om)$, and the nonanticipativity constraints $x=x(\om)$ for $\om \in \Om$ are introduced and then relaxed using Lagrangian multipliers. This allows for scenario-wise decomposition; however, the resulting subproblem in each scenario is a mixed-integer nonconvex program, which can be computationally demanding to solve. 

Recently, smoothing-based approaches have been presented for the continuous case of the min-min problem~\cref{eq:btsp}, where $\X = \R^{n_1}_+$ and $\Y = \R^{n_2}_+$. 
For example, \citet{Liu2020} propose an algorithm that converges almost surely to a generalized critical point. Their approach involves regularizing the recourse function, deriving a \textit{difference-of-convex} decomposition of this regularized function, and then obtaining a convex upper approximation. This allows them to update a solution $x$ by solving approximate problems using the convex approximation. For the convergence analysis, they identify an \textit{implicit convex-concave} property of the recourse function $Q$.
\citet{Li2024} also present a decomposition algorithm for two-stage stochastic programs with implicit convex-concave recourse functions where convex approximations of the recourse function are iteratively generated and solved, and the solutions converge to a critical point.
\citet{Bomze2022} propose a bounding method for a class of bi-parameterized two-stage stochastic nonconvex programs where the objective function involves nonconvex quadratic terms and the feasible region is defined by a simplex.
Note that none of the above approaches is directly applicable to BTSPs with integrality restrictions.

\subsection{Applications of Bi-Parameterized Stochastic Programs} \label{sec:intro-btsp-reform}

Despite the challenges involved in solving BTSPs, they have important applications to various classes of problems. In particular, these include interdiction models and stochastic optimization problems with decision-dependent uncertainty, as described below. Complete details on these BTSP applications are provided in \Cref{apx:btsp-reform}.
\begin{itemize}
    \item \textbf{Bi-parameterized network interdiction models.}
        Network interdiction problems involve sequential adversarial games between two players: an interdictor and a network user. These problems have diverse applications in practice, such as disrupting illicit supply networks~\citep{Morton2007,malaviya2012}, and analyzing vulnerabilities in critical infrastructure~\citep{brown2006}. Stochastic network interdiction problems can be formulated as the following min-max form:
        \begin{equation*}
            \min_{x \in X} \E\Big[ \max_{y \in Y(x, \om)} f(x, y, \om) \Big],
        \end{equation*}
        where $x$ denotes interdictor's decision variables and $y$ represents network user's decision variables. The function $f(x, y, \om)$ and set $Y(x, \om)$ for $\om \in \Om$  are the network user's objective function and feasible set, respectively, and both of them depend on $x$. 
        The formulation of bi-parameterized interdiction models has been considered a general form of interdiction problems in the literature~\citep{SongSmithSurvey}; however, to the best of our knowledge, no standard decomposition approaches have been established to address this general setting.
        
    \item \textbf{Stochastic optimization with decision-dependent uncertainty.} 
        Stochastic optimization often involves decision-dependent uncertainty, where a decision affects the distribution (for example, see \citep{dupacova2006,Hellemo2018}). Decision-dependent two-stage (risk-neutral) stochastic programs can be formulated as:
        \begin{equation*}
            \min_{x \in X}
            \Big\{
                f(x) + \E_{p(x)} \big[Q^s(x, \om) \big]
            \Big\},
        \end{equation*}
        where $Q^s(x, \om) = \min \{ q(\om)^\top y: W(\om) y = r(\om) - T(\om) x,\ y \in \Y \}$ is the recourse function of the classical single-parameterized problem, for each $\om \in \Om$. Here, the expectation is taken with respect to the decision-dependent distribution $p(x) = (p(x, \om))_{\om \in \Om}^\top$, which we denote by $\E_{p(x)}$. By incorporating $p(x, \om)$ into the objective function of each recourse problem $Q^s(x, \om)$, resulting in $(p(x, \om) q(\om)^\top y)$, the above formulation can be easily reformulated into the min-min problem~\cref{eq:btsp}.

        Furthermore, we can also consider a risk measure to model risk-averse behaviors of the decision-maker with the decision-dependent probabilities.
        For instance, decision-dependent two-stage risk-averse stochastic programs with conditional value-at-risk (CVaR) can be formulated as:
        \begin{equation*}
            \min_{x \in X}
            \Big\{
                f(x) + \cvar_\alpha\big(Q^s(x, \om), p(x)\big)
            \Big\},
        \end{equation*}
        where the CVaR is given by (\citet{rockafellar2000}):
        \begin{equation}\label{eq:def-cvar}
        \begin{aligned}
            \cvar_\alpha (Q^s(x, \om), p(x)) = \min_{\eta \in \R} \eta + \frac{1}{1-\alpha}\E_{p(x)}\Big[(Q^s(x, \om) - \eta)_+\Big],
        \end{aligned}
        \end{equation}
        which is a special case of the min-min problem~\cref{eq:btsp}.

    \item \textbf{Distributionally robust optimization with decision-dependent ambiguity.} 
        Distributionally robust optimization (DRO) is an optimization framework that addresses uncertain parameters whose true distribution is unknown. In a DRO model, we consider a set of potential distributions---referred to as \textit{ambiguity set}---and make decisions that hedge against worst-case distributions within this set. The ambiguity set may depend on the first-stage decision $x$; for example, see \citep{Basciftci,Luo2020,Yu,KangBansal2025MathProg}. A DRO model with decision-dependent ambiguity set can be formulated as:
        \begin{equation*}
            \min_{x \in X}
            \Big\{
                f(x) + \max_{p \in \P(x)} \E_p \big[Q^s(x, \om) \big]
            \Big\},
        \end{equation*}
        where $\P(x)$ is a decision-dependent ambiguity set and $Q^s(x, \om)$ is defined as the single-parameterized recourse function for $\om \in \Om$. When the number of scenarios is finite, i.e. $|\Om| < \infty$, the ambiguity set $\P(x)$ can be often represented by a polyhedron (e.g., Wasserstein ambiguity set, moment-matching ambiguity set, and $\phi$-divergence ambiguity set). In this case, by taking dual of the continuous relaxations of recourse problems, we can derive a min-max (BTSP) reformulation of the DRO problem, which is a special case of the min-max problem~\cref{eq:btsp}.

        \paragraph{Non-Relatively Complete Ambiguity Set.} 
        A typical method for decision-dependent DRO problems is to dualize the inner maximization and reformulate the entire problem as a single minimization problem. This dual-based method, however, requires the ambiguity set $\P(x)$ to be nonempty for all $x \in X$. In this paper, we refer to such ambiguity sets as ``relatively complete.'' 
        However, for certain types of ambiguity sets, this nonemptiness assumption can be impractical; for example, when the ambiguity set is constructed to match the moment information of sample data, some first-stage decisions may change the ambiguity set in such a way that moment matching becomes impossible, resulting in an empty ambiguity set (see the computational results in \citet{Yu}). In \Cref{sec:test_result-dddro}, we demonstrate that the dual-based approach fails to solve instances where this occurs. In contrast, the method presented in this paper for BTSPs effectively addresses this issue, providing optimal solutions to DRO problems with ``non-relatively complete'' ambiguity sets.
\end{itemize}

\subsection{Contributions and Organization of this Paper}

In this paper, we introduce BTSPs, where first-stage decisions affect both the objective function and constraints in the recourse problem. To the best of our knowledge, we are the first to study these models with integrality restrictions in both the first and second stages. Our approach can address a wide range of problem classes, allowing $Q$ to be a minimization or a maximization problem, and the first-stage variables $x$ can be binary or mixed-integer. Specifically, we present two exact algorithms for the min-max and min-min cases when $x$ is binary, and a regularization-augmented algorithm for the mixed-binary case. The numerical results demonstrate the efficiency of our methods. They solve all tested instances of a stochastic network interdiction problem in 23 seconds on average, while a benchmark method fails to solve any instances within the time limit of 3600 seconds. Additionally, our method outperforms existing methods for a stochastic facility location problem, achieving optimal solutions, on average, 18.4 times faster.
Furthermore, we demonstrate that our approach can effectively address DRO problems with decision-dependent ambiguity sets that are not necessarily relatively complete, i.e., $\P(x)$ may be empty for some $x \in X$; we note that this setting remains unaddressed in the literature~(\cite{Basciftci,Luo2020,Yu}). For these problem instances, our approach achieves optimal solutions, on average, 5.3 times faster than an existing approach.

The remainder of the paper is organized as follows. In \Cref{sec:method_min-max}, we introduce an exact method for the min-max problem. In \Cref{sec:method_min-min}, we present its extension for solving the min-min problem. In \Cref{sec:method_reg}, we propose a regularization-augmented variant that addresses BTSPs with mixed-integer first-stage variables. Computational results are provided in \Cref{sec:test_result}, and concluding remarks are given in \Cref{sec:conclusion}. For convenience, all proofs are presented in \Cref{apx:proofs}.

We make the following assumption on the feasible sets throughout the paper:
\begin{assumption} \label{assumption:complete_recourse}
    The set $X$ is nonempty and compact. Also, for all $x \in X$, the sets $Y(x, \om)$ for $\om \in \Om$ are nonempty and compact.
\end{assumption}

\paragraph{Notation:} Let $[n] := \{1, 2, \dots, n\}$ for a positive integer $n$.

\section{An Exact Algorithm for the Min-Max BTSP} \label{sec:method_min-max}

In this section, we introduce a decomposition method referred to as the Lagrangian-integrated $L$-shaped ($L^2$) method to exactly solve the min-max BTSP~\cref{eq:btsp} with integer variables in both stages. 

\subsection{A Lagrangian Dual Formulation}

We begin by deriving a Lagrangian dual of the min-max problem.
To do so, we first rewrite the recourse problem~\cref{eq:Q} in the following form by introducing proxy variables $z$ for $x$:
\begin{subequations} \label{eq:Q_om_re}
\begin{align}
    Q(x, \om) = \max \ & q(\om)^\top y + x^\top G(\om) y  \\
    \text{s.t.} \ 
    & T(\om) z + W(\om) y = r(\om), \ z \in X, \ y \in \Y,\\
    & z = x. \label{eq:aux_con}
\end{align}
\end{subequations}
Note that we substitute $x$ with $z$ only in the constraints, while leaving $x$ unchanged in the objective---ensuring that the problem remains linear for fixed $x$.
We then relax constraints~\cref{eq:aux_con} with Lagrangian multipliers $\lambda(\om) \in \R^{n_1}$. The Lagrangian relaxation is given by $(L(x, \lambda(\om), \om) - \lambda(\om)^\top x)$ where
\begin{subequations} \label{eq:lag_Q_om}
\begin{align}
    L(x, \lambda(\om), \om) =
    \max \ & q(\om)^\top y + x^\top G(\om) y + \lambda(\om)^\top z \\
    \text{s.t.} \ 
    & T(\om) z + W(\om) y = r(\om), z \in X, \ y \in \Y. \label{eq:lag_Q_om_c}
\end{align}
\end{subequations}
For any $\lambda(\om)\in\R^{n_1}$, the Lagrangian relaxation $(L(x, \lambda(\om), \om) - \lambda(\om)^\top x)$ provides an upper bound on $Q(x, \om)$. Thus, the tightest upper bound is obtained by solving the following Lagrangian dual:
\begin{equation} \label{eq:dual_Q_om}
    D(x, \om) := \min_{\lambda(\om) \in \R^{n_1}} \Big\{ L(x, \lambda(\om), \om) - \lambda(\om)^\top x \Big\}.
\end{equation}

\begin{assumption} \label{assumption:binary_x_component_minmax}
    In the min-max problem~\cref{eq:btsp}, variables $x$ are mixed-integer. Let $I \subseteq [n_1]$ denote the index set corresponding to the components of $x$ that affect the recourse function $Q(x, \om)$. We assume that for all $i \in I$, $x_i$ is binary.
\end{assumption}

Under the above assumption, we show that strong duality holds for the Lagrangian dual. It is worth noting that, although similar duality results have been established in the literature (e.g., \cite{zou2019}), our result is distinct in that we consider bi-parameterized problems and relax only the constraint parameterization, while leaving $x$ in the objective function.

\begin{theorem} \label{thm:strong_duality_Q_om}
    Under \cref{assumption:binary_x_component_minmax}, $Q(x, \om) = D(x, \om)$ for any $\om\in\Om$ and $x \in X$, i.e., strong duality holds for the Lagrangian dual~\cref{eq:dual_Q_om}.
\end{theorem}

Consequently, the min-max problem~\cref{eq:btsp} can be exactly reformulated as
\begin{equation} \label{eq:min-max_bilinear}
    \min \bigg\{
        f(x) + \sum_{\om \in \Om} p(\om) \Big( L(x, \lambda(\om), \om) - \lambda(\om)^\top x \Big)
        : x \in X, \ \lambda(\om) \in \R^{n_1}, \ \forall \om \in \Om
    \bigg\}.
\end{equation}
An advantage of addressing the reformulated problem instead of the original problem is that the Lagrangian function $L(\cdot, \cdot, \om)$ for each $\om \in \Om$ is jointly convex on $\conv(X) \times \R^{n_1}$. This convexity enables the use of a subgradient method, which decomposes the computation of $L$. However, nonconvex terms $\lambda(\om)^\top x$ appear in the objective function, which can cause computational challenges. 

To address this issue, we introduce an analytical form for the optimal Lagrangian multipliers.
\begin{lemma} \label{lem:sigma_optimal_condition}

There exist real values $\sigma_i \in \R_+$ for $i \in I$ such that, for any $\hat{x} \in X$ and $\om \in \Om$, the solution $\hat{\lambda}(\om) = (\hat{\lambda}_{1}(\om), \dots, \hat{\lambda}_{n_1}(\om))$, defined by
\begin{equation}
    \hat{\lambda}_{i}(\om) = \begin{cases}
        -\sigma_i & \text{ if } i \in I \text{ and } \hat{x}_i = 0 \\
        \sigma_i & \text{ if } i \in I \text{ and } \hat{x}_i = 1 \\
        0 & \text{ otherwise}
    \end{cases},
\end{equation} 
is optimal to the Lagrangian dual~\cref{eq:dual_Q_om}. Here, $I \subseteq [n_1]$ is the index set corresponding to the components of $x$ that affect the recourse function $Q(x, \om)$.

Furthermore, the following provides a sufficient condition on the values of $\sigma$ to guarantee the above optimality:
\begin{equation} \label{eq:sigma_condition}
    \sum_{i \in I: \hat{z}_i (\om) \neq \hat{x}_i} \sigma_i \geq 
    \max_{y \in Y(\hat{x}, \om)} 
    \Big\{
    (q(\om) + G(\om)^\top \hat{x})^\top (\hat{y}(\om) - y) 
    \Big\},
\end{equation}
where $(\hat{y}(\om), \hat{z}(\om))$ is an optimal solution of the Lagrangian relaxation~\cref{eq:lag_Q_om} given $(\hat{x}, \hat{\lambda}(\om), \om)$.
\end{lemma}

Using this analytical form for an optimal solution, we can derive an exact reformulation of the min-max problem as follows.

\begin{theorem} \label{thm:min-max_reform}
    Let $(\sigma_i)_{i \in I}$ satisfy the condition~\cref{eq:sigma_condition} for all $\hat{x} \in X$ and $\om \in \Om$.
    The following problem is an exact reformulation of the min-max problem~\cref{eq:btsp}:
    \begin{equation} \label{eq:min-max_reform}
        \min_{x \in X} 
        \Big\{
        f(x) - \sum_{i \in I} \sigma_i x_i +
            \sum_{\om \in \Om} p(\om) \hat{Q}(x, \om)
        \Big\}
    \end{equation}
    where, for each $\om \in \Om$,
    \begin{equation} \label{eq:min-max_reform_Q}
    \begin{aligned}
        \hat{Q}(x, \om) = 
        \max \ & q(\om)^\top y + x^\top G(\om) y + \sum_{i \in I} \sigma_i (2 x_i - 1) z_{i} \\
        \text{s.t.} \ 
        & T(\om) z + W(\om) y = r(\om), \ z \in X, \ y \in \Y.
    \end{aligned}
    \end{equation}
\end{theorem}

\begin{remark}

To find values $\sigma_i, i \in I,$ that satisfy condition~\cref{eq:sigma_condition} for all $\hat{x} \in X$ and $\om \in \Om$, we can use the following upper and lower bounds: $\bar{UB} \geq \max_{y \in Y(x, \om)} (q(\om) + G(\om)^\top x)^\top y$ and $\bar{LB} \leq \min_{y \in Y(x, \om)} (q(\om) + G(\om)^\top x)^\top y$, respectively. Notice that their difference is always greater than or equal to the right-hand side of the condition \cref{eq:sigma_condition}, that is,
$
    \max_{y \in Y(\hat{x}, \om)} (q(\om) + G(\om)^\top \hat{x})^\top (\hat{y}(\om) - y)
    = (q(\om) + G(\om)^\top \hat{x})^\top\hat{y}(\om) - \min_{y \in Y(\hat{x}, \om)} (q(\om) + G(\om)^\top \hat{x})^\top y 
    \leq \bar{UB} - \bar{LB}.
$
By setting $\sigma_i$ equal to this difference $(\bar{UB} - \bar{LB})$ for $i \in I$, we can ensure that the condition~\cref{eq:sigma_condition} is satisfied. 

Alternatively, we may use an optimization-based approach to compute these values. For example, define $\Delta_y(x,\om) = \{y^1 - y^2: y^1 , y^2 \in Y(x, \om)\}$ and consider
$
    \sigma_i(\om) = \max 
    \big\{
        \big(\frac{q(\om)}{|I|} + G_i(\om) x_i\big)^\top \delta_y :
        x \in X,\ 
        \delta_y \in \Delta_y(x, \om)
    \big\}
$
for $i \in I$ and $\om \in \Om$, where $G_i(\om)$ denotes the $i$th column of $G(\om)$. Notice that
$
    \sum_{i \in I} \sigma_i(\om) 
    \geq
    \max
    \{
        (q(\om) + G(\om)^\top x)^\top \delta_y :
        x \in X,\ 
        \delta_y \in \Delta_y(x, \om)
    \}
$
for $\om \in \Om$. This inequality holds since, for an optimal solution $(x^\ast, \delta_y^\ast)$ to the right-hand side problem, its component $(x^\ast_i, \delta^\ast_y)$ is feasible to the corresponding $i$th problem defining $\sigma_i(\om)$. These $\sigma_i(\om)$ values ensure that the condition~\cref{eq:sigma_condition} is satisfied for each $\om \in \Om$. By taking the maximum over $\om \in \Om$, i.e., $\sigma_i = \max_{\om \in \Om} \sigma_i(\om)$, we have $\sigma_i$ values that satisfy the condition~\cref{eq:sigma_condition} for all $\hat{x} \in X$ and $\om \in \Om$.
$\blacksquare$\end{remark}

\begin{algorithm}[t]
\caption{The $L^2$ method for the min-max problem}\label{alg:min-max-L2}
\textbf{Initialization}: $l \gets 1, LB \gets -\infty, UB \gets \infty, x^1 \in X$\;
\While{$(UB - LB) / UB > \epsilon^{tol}$}{
    Solve the subproblems, obtain $\{(y^l(\om), z^l(\om))\}_{\om \in \Om}$, and compute $\hat{Q}(x^{l}, \om)$\; \label{line:min-max-subproblem}
    Compute $\phi^l = f(x^{l}) - \sum_{i \in I} \sigma_i x^{l}_i + \sum_{\om \in \Om} p(\om) \hat{Q}(x^{l}, \om)$\;
    \If{$\phi^l < UB$}{
        $UB \gets \phi^l$; $x^\ast \gets x^{l}$\;
    }
    % \textbf{if} $\phi^l < UB$ \ \textbf{then} \  $UB \gets \phi^l$; $x^\ast \gets x^{l}$\;
    Add an optimality cut to master problem\;
    Solve the master problem and obtain optimal solution $(x^{l+1}, \hat{\theta})$\;
    $LB \gets f(x^{l+1}) - \sum_{i \in I} \sigma_i x^{l+1}_i + \hat{\theta}$; \ \label{line:min-max_compute_lower_bound} $l \gets l + 1$\;
}    
\Return{$x^\ast$}\;
\end{algorithm}

\subsection{The Lagrangian-Integrated \texorpdfstring{$L$}{L}-Shaped Method}

We now present the $L^2$ method for solving the min-max problem~\cref{eq:btsp} through its reformulation~\cref{eq:min-max_reform}. At each iteration, the $L^2$ method adds a valid optimality cut that approximates the function $\sum_{\om \in \Om} p(\om) \hat{Q}(x, \om)$ in \cref{eq:min-max_reform}.
The pseudo-code of the $L^2$ method is outlined in \Cref{alg:min-max-L2}. It starts by initializing iteration counter $l$ to $1$, lower bound $LB$ to $-\infty$, and upper bound $UB$ to $\infty$. Let $x^1 \in X$ be any feasible solution. For each iteration $l$, we first solve the subproblems for $\om \in \Om$, which is given by the formulation~\cref{eq:min-max_reform_Q} with $x = x^l$.
Subsequently, we compute the objective value at the current solution $x^{l}$, denoted by $\phi^l$. If this value $\phi^l$ is less than the current upper bound $UB$, then it replaces $UB$, and the current solution is saved as the best-known solution. 
Next, using the optimal solutions to the subproblems, an optimality cut is obtained and then added to the master problem. In iteration $l$, master problem is given as follows: 
\begin{equation}\label{eq:min-max_master}
\begin{aligned}
    \min \ \bigg\{ & f(x) - \sum_{i \in I} \sigma_i x_i + \theta: x \in X, \\
        & 
        \theta \geq \sum_{\om \in \Om} p(\om) \Big(
            q(\om)^\top y^k(\om) + x^\top G(\om) y^k(\om)
            + \sum_{i \in [n_1]} \sigma_i  (2x_i - 1) z^k_{i}(\om)
        \Big), \ \forall k \in [l] \bigg\},
\end{aligned}
\end{equation}
where $(y^k(\om), z^k(\om))_{k \in [l]}$ are the optimal solutions of the subproblems at the corresponding iterations.
In the following step, we solve the master problem to obtain a first-stage solution $x^{l+1}$ and the updated lower bound $LB$. We repeat this procedure until the optimality gap is less than or equal to a predetermined tolerance level $\epsilon^{tol}$.

\begin{proposition}\label{prop:finite_convergence}
    Let values $\sigma_i \in \R_+, i \in I,$ satisfy the condition~\cref{eq:sigma_condition} for all $x \in X$ and $\om \in \Om$, and let $\epsilon^{tol} = 0$.
    Under \Cref{assumption:binary_x_component_minmax}, \Cref{alg:min-max-L2} terminates in a finite number of iterations with $x^\ast$ optimal to the min-max BTSP~\cref{eq:btsp}.
\end{proposition}

\begin{remark} \label{rem:L2-fixed-z}
    To accelerate \Cref{alg:min-max-L2}, we can fix $z$ to $x^l$ when solving the subproblem in \Cref{line:min-max-subproblem}. Solving this restricted subproblem with $z = x^l$ yields a feasible solution to the original problem~\cref{eq:min-max_reform_Q}, allowing us to derive an optimality cut valid for $\hat{Q}(x, \om)$. This approach reduces computational effort while still providing a valid cut.
$\blacksquare$\end{remark}

\section{An Exact Algorithm for the Min-Min BTSP} \label{sec:method_min-min}

We now present the $L^2$ method applied to the min-min problem~\cref{eq:btsp}. To enable the application of the $L^2$ method described in \Cref{sec:method_min-max}, we convexify the recourse problem by incorporating parametric inequalities---as established in, e.g., \citet{gade2014}---and derive a min-max reformulation. We show that, with this convexification of the bi-parameterized recourse, optimality cuts remain valid.

For the min-min problems, we assume that $x$ is binary to ensure the exactness of our algorithm.
\begin{assumption} \label{assumption:binary_x_component_minmin}
    In the min-min problem, all first-stage variables $x$ are binary, i.e., $\X = \{0,1\}^{n_1}.$
\end{assumption}

Let $S(\om) := \{ (x, y) \in X \times \Y: T(\om) x + W(\om) y = r(\om) \}$ be a lifted feasible set in the $(x, y)$-space for $\om \in \Om$. Assume that $\tilde{T}(\om) \in \Qu^{m(\om) \times n_1}, \tilde{W}(\om) \in \Qu^{m(\om) \times n_2},$ and $\tilde{r}(\om) \in \Qu^{m(\om)}$ are appropriately sized matrices and vectors such that $\conv(S(\om)) = \{(x, y) \in \R^{n_1}_+ \times \R^{n_2}_+: \tilde{T}(\om) x + \tilde{W}(\om) y \geq \tilde{r}(\om) \}$. 
The intersection of $\conv(S(\om))$ and the hyperplane $\{(x, y): x = \hat{x}\}$ for any $\hat{x} \in X$ is a face of $S(\om)$, as $x$ is binary, thereby all extreme points in the intersection have the $y$ components in $\Y$.
By this observation, for any $\hat{x} \in X$, we have 
\begin{align}
     Q(\hat{x}, \om) 
     &= \min \bigg\{ 
         q(\om)^\top y + \hat{x}^\top G(\om) y : (x, y) \in \conv(S(\om)) \cap \Big\{(x, y) : x = \hat{x}\Big\}
    \bigg\} \\
    &= \min \bigg\{ 
        q(\om)^\top y + \hat{x}^\top G(\om) y : \tilde{W}(\om) y \geq \tilde{r}(\om) - \tilde{T}(\om) \hat{x},\ 
        y \in \R^{n_2}_+
    \bigg\}.\label{eq:min-min_Q_conv}
\end{align} 
Since $Y(x, \om)$ is bounded, we can assume, without loss of generality, that the number of constraints $m(\om)$ in \cref{eq:min-min_Q_conv} is finite. That is, the convex hull $\conv(S(\om))$ can be represented with finitely many parametric inequalities in the form of $(\alpha^1)^\top x + (\alpha^2)^\top y \geq \beta$; e.g., Lift-and-Project cuts and Gomory cuts.  
Thus, the linear programming dual of \cref{eq:min-min_Q_conv} is given by
\begin{equation} \label{eq:min-min_Q_dual}       
    Q(x, \om) = \max \bigg\{ 
        \pi^\top (\tilde{r}(\om) - \tilde{T}(\om) x) :
        \tilde{W}(\om)^\top \pi \leq q(\om) + G(\om)^\top x, \ \pi \in \R_+^{m(\om)}
    \bigg\}.
\end{equation}
Then, we can rewrite the min-min problem as the following min-max problem:
\begin{equation} \label{eq:min-min_dual}
    \min_{x \in X} \bigg\{
        f(x) + \sum_{\om \in \Om} p(\om) \max \Big\{
            \pi^\top (\tilde{r}(\om) - \tilde{T}(\om) x) : 
            \tilde{W}(\om)^\top \pi \leq q(\om) + G(\om)^\top x, \ \pi \in \R_+^{m(\om)}
        \Big\}
    \bigg\}.
\end{equation}
By applying the result in \Cref{thm:min-max_reform} to the above min-max reformulation, we obtain the following another reformulation of the min-min problem: $\min_{x \in X} f(x) - \sigma^\top x + \sum_{\om \in \Om} p(\om) \hat{Q}(x, \om)$
where 
\begin{equation} \label{eq:min-min_reform_Q}
\begin{aligned}
    \hat{Q}(x, \om) = 
    \max \bigg\{ & \pi^\top (\tilde{r}(\om) - \tilde{T}(\om) x) + \sum_{i \in [n_1]} \sigma_i (2x_i - 1) z_{i} : \\
    & \tilde{W}(\om)^\top \pi - G(\om)^\top z \leq q(\om), \ z \in X, \ \pi \in \R^{m(\om)}_+ \bigg\}.
\end{aligned}
\end{equation} 

We now present the $L^2$ method for the min-min problem. Its pseudo-code is outlined in \Cref{alg:min-min-L2}. The algorithm sequentially convexifies the set $S(\om)$, i.e., it constructs the information $\{(\tilde{T}(\om), \tilde{W}(\om), \tilde{r}(\om))\}_{\om \in \Om}$ by adding parametric inequalities to subproblems. To this end, we consider an oracle that provides a violated parametric inequality if there exists any. We refer to such an oracle as the \textit{cut-generating} oracle.

\begin{algorithm}[t]
\caption{The $L^2$ method for the min-min problem}\label{alg:min-min-L2}
\textbf{Initialization}: $l \gets 1, LB \gets -\infty, UB \gets \infty, x^1 \in X$\;
\While{$(UB - LB) / UB > \epsilon^{tol}$}{
    \For{$\om \in \Om$}{
        Solve the primal subproblem and obtain optimal solution $y^l(\om)$\; \label{line:min-min_solve_primal_sub}
        \If{$y^l(\om) \notin \Y$}{
            Find a violated parametric inequality using the cut-generating oracle\;
            Augment $(T^{l-1}(\om), W^{l-1}(\om), r^{l-1}(\om))$ with the data of the violated inequality\;
        }
        $(T^l(\om), W^l(\om), r^l(\om)) \gets (T^{l-1}(\om), W^{l-1}(\om), r^{l-1}(\om))$\;
    }
    Solve the dual subproblems, obtain $\{(\pi^l(\om), z^l(\om))\}_{\om \in \Om}$, and compute $\hat{Q}^l(x^{l}, \om)$\; \label{line:min-min_solve_dual_sub}
    
    Compute $\psi^l = f(x^{l}) - \sigma^\top x^{l} + \sum_{\om \in \Om} p(\om) \hat{Q}^l(x^{l}, \om)$\;
    \If{$y^l(\om) \in \Y, \forall \om \in \Om$ and $\psi^l < UB$}{
        $UB \gets \psi^l$; $x^\ast \gets x^{l}$\;
    }
    % \textbf{if} $y^l(\om) \in \Y, \forall \om \in \Om$ and $\psi^l < UB$ \ \textbf{then} \ 
        $UB \gets \psi^l$; $x^\ast \gets x^{l}$\;
    Add an optimality cut to master problem\;
    Solve the master problem and obtain optimal solution $(x^{l+1}, \hat{\theta})$\;
    $LB \gets f(x^{l+1}) - \sigma^\top x^{l+1} + \hat{\theta}$; \ $l \gets l + 1$\;
}    
\Return{$x^\ast$}\;
\end{algorithm}

\Cref{alg:min-min-L2} starts by initializing iteration counter $l$ to 1, lower bound $LB$ to $-\infty$, and upper bound $UB$ to $\infty$. Let $x^1 \in X$ be any feasible solution. At each iteration $l$, we first solve the following problems, called \textit{primal subproblems}:
\begin{equation} \label{eq:min-min_primal_sub}
    \min \Big\{
        q(\om)^\top y + (x^{l})^\top G(\om) y : W^{l-1}(\om) y \geq r^{l-1}(\om) - T^{l-1}(\om) x^{l}, \ 
        y \in \R^{n_2}_+
    \Big\},
\end{equation}
where $(T^0(\om), W^0(\om), r^0(\om)) = (T(\om), W(\om), r(\om))$, for $\om \in \Om$.
The matrices and vector are updated to $(T^{l}(\om), W^{l}(\om), r^{l}(\om))$ through the following procedure. Let $y^l(\om)$ be the optimal solution to the primal subproblems for $\om \in \Om$. If $y^l(\om) \notin \Y$, then the cut-generating oracle is used, and a violated parametric inequality is added to the primal subproblem. Otherwise, we keep the current matrices and vector. Notice that the primal subproblem is a relaxation of the recourse problem~\cref{eq:min-min_Q_conv}, thereby through this approach the $L^2$ method refines the relaxation iteratively. In the next step, we solve the following problems for $\om \in \Om$, referred to as \textit{dual subproblems}:
\begin{equation} \label{eq:min-min_dual_sub}
\begin{aligned}
    \hat{Q}^l (x^{l}, \om) = 
    \max \Big\{ 
        &
        \pi^\top(r^l(\om) - T^l(\om) x^{l}) + \sum_{i \in [n_1]}\sigma_i (2x^{l}_i - 1) z_{i} : \\
        & W^l(\om)^\top \pi - G(\om)^\top z \leq q(\om),\ 
        z \in X, \ \pi \in \R^{m^l}_+
    \Big\}.
\end{aligned}
\end{equation}
where $\sigma$ is a given vector that satisfies the condition~\cref{eq:sigma_condition} for all $\om \in \Om$ and $x \in X$.
Note that the dimension $m^l$ of variables $\pi$ varies over iterations corresponding to the updates in the information $(T^l(\om), W^l(\om), r^l(\om))$. 
After solving the dual subproblems, an under-approximation $\psi^l$ of the objective value at the current solution $x^{l}$ is computed; this objective value is exact if $y^l(\om) \in \Y$ for all $\om \in \Om$. If $\psi^l$ is exact and $\psi^l < UB$, then we update the upper bound, and the current solution is marked as the best solution so far. Subsequently, using the optimal solutions of the dual subproblems, an optimality cut is computed and added to \textit{master} problem that is given by 
\begin{subequations} \label{eq:min-min_master}
\begin{align}
    \min \ & \bigg\{f(x) - \sigma^\top x + \theta: x \in X, \\
        & 
        \begin{aligned}
        \theta \geq & \sum_{\om \in \Om} p(\om) \Big(
            \pi^{k}(\om)^\top (r^k(\om) - T^k(\om) x) + \sum_{i \in [n_1]} \sigma_i z_{i}^k(\om) (2 x_i - 1)
        \Big), \quad \forall k \in [l]\bigg\}. \label{eq:min-min_cut}
        \end{aligned}
\end{align}
\end{subequations} 
Constraints~\cref{eq:min-min_cut} are optimality cuts added for iteration $k \in [l]$ where $(\pi^k(\om), z^k(\om))_{k \in [l]}$ are optimal solutions of the dual subproblems at the corresponding iterations.
\begin{proposition} \label{prop:min-min_valid_cut}
    The optimality cut obtained at iteration $l$ is valid for $\sum_{\om\in\Om} p(\om) \hat{Q}(x, \om)$, i.e., $\sum_{\om\in\Om} p(\om) \hat{Q}(x, \om) \geq \sum_{\om \in \Om} p(\om) \Big(
            \pi^l(\om)^\top (r^l(\om) - T^l(\om) x) + \sum_{i \in [n_1]} \sigma_i z_{i}^l(\om) (2 x_i - 1)  
        \Big)$.
\end{proposition}

Next, we solve the master problem and obtain a new lower bound on the optimal objective value. It also identifies a solution $x^{l+1}$ to be explored in the next iteration. This process is repeated until the optimality gap, $(UB-LB)/UB$, equals or falls below a predetermined tolerance level $\epsilon^{tol} \geq 0$.

\section{A Regularization-Augmented Algorithm for BTSPs} \label{sec:method_reg}

The proposed exact methods (\Cref{alg:min-max-L2,alg:min-min-L2}) require a predetermined value of $\sigma$ that satisfies the condition~\cref{eq:sigma_condition}. In this section, we introduce an alternative approach that does not require specifying $\sigma$ a priori. Instead of relying on the analytical form in \Cref{lem:sigma_optimal_condition}, we directly address the bilinear problem~\cref{eq:min-max_bilinear}. Specifically, we propose a subgradient method and improve its computational efficiency by stabilizing updates of the Lagrangian multipliers through reformulation and regularization techniques. This method provides an approximate solution for BTSPs with general mixed-binary $x$, that is, for the cases where \Cref{assumption:binary_x_component_minmax,assumption:binary_x_component_minmin} do not hold.

For simplicity, in this section, we focus on the min-max problem where first-stage variables $x$ are general mixed-binary and the functions $L(x, \lambda(\om), \om)$ for $\om \in \Om$ are defined as \cref{eq:lag_Q_om}. We note that the proposed method can be adapted to the min-min problem in a similar manner after taking the min-max reformulation, described in \Cref{sec:method_min-min}.

We begin by discussing how stabilizing multiplier updates can improve computational efficiency.
Consider a standard subgradient method that iteratively updates both $x$ and $\lambda(\om)$. An issue arises in updating $\lambda(\om)$, depending on the values of $x$. 
For example, suppose $x$ is pure binary and fixed at $\hat{x}$, and that the subgradient method has converged to the optimal multipliers $\lambda_i(\om) = \sigma_i (2\hat{x}_i - 1)$ for $i \in [n_1]$, where $\sigma$ satisfies condition~\cref{eq:sigma_condition}. Now, consider $i \in [n_1]$ such that $\hat{x}_i = 0$, so $\lambda_i(\om) = -\sigma_i$. If, in the next iteration, $\hat{x}_i$ is switched to $1$, the optimal multiplier changes to $\sigma_i$. This requires a change of $2\sigma_i$ in the multiplier, which may demand many update steps and lead to slow convergence.

Motivated by this observation, we propose a reformulation that prevents sign flips in updating targets. The key idea is to set the multipliers to be nonpositive when $x_i = 0$ and nonnegative when $x_i = 1$ for the binary variables and to update their magnitudes at each iteration, rather than updating the multipliers directly. Let $N' \subseteq [n_1]$ be the index set of binary variables in $x$, and let $N'' = [n_1]\setminus N'$. The reformulation is given by
\begin{subequations} \label{eq:min-max_bilinear-reform}
\begin{align}
    \min_{x \in X} \ &
        f(x) + \sum_{\om \in \Om} p(\om) \Big( L(x, \lambda(\om), \om) - \lambda(\om)^\top x \Big) \\
    \text{s.t.} \ 
    & \lambda_{i}(\om) = \mu_{i}(\om) (2x_i - 1),\ 
    \mu_i(\om) \in \mathbb{R}_+, \ \forall \om \in \Om, i \in N', \ 
    \lambda_i(\om) \in \mathbb{R}, \ \forall \om \in \Om, i \in [n_1]. \label{eq:min-max_bilinear-reform-constr}
\end{align}
\end{subequations}
Here, a variable $\mu_i(\om)$ represents the magnitude of the corresponding multiplier $\lambda_i(\om)$, for each $i \in N'$. Thus, updating this $\mu_i (\om)$ for $i \in N'$ instead of $\lambda_i(\om)$ in a subgradient method can stabilize the update process.
When $x$ is pure binary, the reformulation~\cref{eq:min-max_bilinear-reform} is exact, as an optimal solution, given in \Cref{lem:sigma_optimal_condition}, is feasible to this problem. If $x$ is mixed-binary, then it provides an approximate solution to the original problem. 

\begin{algorithm}[t]
\caption{The regularized $L^2$ method for the min-max problem}\label{alg:min-max_alt}
\textbf{Initialization}: $l \gets 1, LB \gets -\infty, UB \gets \infty, x^1 \in X, \lambda^1(\om) \in \R^{n_1}_+$ for $\om \in \Om$\;
\While{$(UB - LB) / UB > \epsilon^{tol}$}{
    Solve the subproblems~\cref{eq:min-max_sub_alt}, obtain $\{(y^l(\om), z^l(\om))\}_{\om \in \Om}$, and compute $\hat{Q}^{l}(x^{l}, \lambda^{l}(\om), \om)$\; \label{line:min-max_solve_dual_sub_alt}
    Compute $\phi^l = f(x^{l}) + \sum_{\om \in \Om} p(\om) \Big(\hat{Q}^l(x^{l}, \lambda^{l}(\om), \om) - \lambda^{l}(\om)^\top x^{l}\Big)$\;
    \If{$\phi^l < UB$}{
        $UB \gets \phi^l$; $x^\ast \gets x^{l}$\;
    }
    % \textbf{if} $\phi^l < UB$ \ \textbf{then} \ $UB \gets \phi^l$; $x^\ast \gets x^{l}$\;
    Add an optimality cut to master problem~\cref{eq:min-max_master_alt}\;
    Solve master problem~\cref{eq:min-max_master_alt} and obtain $x^{l+1}$, $(\lambda^{l+1}(\om))_{\om \in \Om}$, $(\mu^{l+1}(\om))_{\om \in \Om}$, and $\hat{\theta}$\;
    $LB \gets f(x^{l+1}) + \hat{\theta} - \sum_{\om \in \Om} p(\om) \lambda^{l+1}(\om)^\top x^{l+1} + \gamma R(\mu^{l+1})$; \label{line:compute_lower_bound_alt} $l \gets l + 1$\;
}    
\Return{$x^\ast$}\;
\end{algorithm}

We now present a variant of the $L^2$ method tailored to the above reformulation. We refer to this algorithm as the \textit{regularized $L^2$ method}, as it incorporates a regularization term into the master problem. This regularization term is introduced to make the optimal multiplier unique and prevent having excessively large multiplier values. The pseudo-code of the regularized $L^2$ method is outlined in \Cref{alg:min-max_alt}. We describe its details below for the completeness of the paper. It starts by initializing iteration counter $l$ to $1$, bounds $LB$ to $-\infty$, and $UB$ to $\infty$. Initial feasible solutions are denoted by $x^1$ and $\lambda^1(\om)$ for $\om \in \Om$. 
Next, we solve the subproblems, which are defined as follows, given $x^{l}$ and $\lambda^l(\om)$ for $\om\in\Om$:
\begin{equation} \label{eq:min-max_sub_alt}
\begin{aligned} 
    \hat{Q}^l (x^{l}, \lambda^{l}(\om), \om) = 
    \max \Big\{ 
        &q(\om)^\top y + (x^{l})^\top G(\om) y +  \lambda^{l}(\om)^\top z : \\
        & T(\om) z + W(\om) y = r(\om),\ 
        z \in X, \ y \in \Y
    \Big\}.
\end{aligned}
\end{equation}
Let $(y^l(\om), z^l(\om))$ be the optimal solution of the subproblem for $\om \in \Om$.
Using the optimal objective values of the subproblems, we can obtain an under-approximation $\phi^l$ of the optimal objective value of the min-max problem. This under-approximation can then be used to update the upper bound $UB$.
In the following step, we add an optimality cut to the master problem given as follows:
\begin{subequations} \label{eq:min-max_master_alt}
\begin{align}
    \min_{x \in X} \ & f(x) + \theta - \sum_{\om \in \Om} p(\om) \lambda(\om)^\top x + \gamma 
    R(\mu)
    \label{eq:min-min_master_alt_obj} \\
        \text{s.t.}\ 
        & \cref{eq:min-max_bilinear-reform-constr}, \ 
        \theta \geq \sum_{\om \in \Om} p(\om) \Big(
            q(\om)^\top y^k(\om) + x^\top G(\om) y^k(\om) 
            + \lambda(\om)^\top z^k(\om)
        \Big), \  \forall k \in [l],
\end{align}
\end{subequations}
where $\mu = (\mu(\om))_{\om \in \Om}$. The term $R(\mu)$ denotes a regularization term, such as $\ell_1$ regularization $\norm{\mu}_1$ or Tikhonov regularization $\norm{\mu}_2^2$, with $\norm{\cdot}_p$ denoting the $\ell^p$ norm. Here, $\gamma$ is a predetermined parameter that controls the significance of the regularization term.
Let $x^{l+1}$, $(\lambda^{l+1}(\om))_{\om \in \Om}$, $\mu^{l+1} = (\mu^{l+1}(\om))_{\om \in \Om}$, and $\hat{\theta}$ be an optimal solution of the master problem at iteration $l$.
In \Cref{line:compute_lower_bound_alt}, the lower bound $LB$ is updated using the optimal objective value of the master problem.

\section{Computational Results} \label{sec:test_result}

We conduct numerical experiments to evaluate the computational efficiency of the proposed approaches. We consider three problem sets: \textit{(a)} bi-parameterized min-min stochastic facility location problem, \textit{(b)} bi-parameterized min-max stochastic network interdiction problem, and \textit{(c)} distributionally robust facility location problem with decision-dependent ambiguity set. In our implementation of the $L^2$ methods, $\sigma$ is set to a sufficiently large value, determined through preliminary experiments to ensure that the test instances are solved to optimality. For the regularized $L^2$ methods, we simplify the model by removing the dependency of the variables $\mu(\om)$ on $\om$. Since the problem sets involve pure binary $x$, this modification still ensures exact solutions while reducing both the dimensionality of the solution space and the overall computational burden. Detailed parameter settings are provided for each specific problem in the ensuing sections.
All algorithms were coded in Julia 1.9 and implemented through the branch-and-cut framework of Gurobi 9.5. The optimality tolerance is set to $10^{-4}$, and the time limit is set to 3600 seconds. We conducted all tests on a machine with an Intel Core i7 processor (3.8 GHz) and 32 GB of RAM, using a single thread.

\subsection{Bi-Parameterized (Min-Min) Facility Location Problem}

We introduce a bi-parameterized facility location problem (BiFLP), where the first-stage decision involves both establishing facilities and contracting outsourcing suppliers. Unlike the traditional facility location problem, customer demand can also be met by outsourcing suppliers, with contracts established in advance to reduce procurement costs. The decision-maker must balance the trade-off between building their own facilities (which incur higher fixed costs but lower variable costs) and outsourcing (which involves lower fixed costs but higher variable costs).

Let $I=\{1, \dots, n_1\}$ be the set of potential facility locations, $J=\{1, \dots, n_2\}$ the set of demand locations, and $K=\{1, \dots, n_3\}$ the set of potential outsourcing suppliers. Binary variables $x'_i$ for $i \in I$ and $x''_k$ for $k \in K$ represent the decision to build a facility at location $i$ and the decision to establish an outsourcing contract with supplier $k$, respectively, subject to budget $b > 0$. Let $x = (x', x'')$ denote a vector of all first-stage decision variables. The random demand at location $j \in J$ is represented by random variable $\tilde{d}_j$, and its realizations are denoted by $d_{j}(\om)$ for $\om \in \Om$. Demand can be fulfilled by both facility at $i \in I$ and supplier $k \in K$. In scenario $\om \in \Om$, the flow from facility $i$ to demand location $j$ is denoted by variable $y_{ij}(\om)$, and from supplier $k$ to $j$ by variable $u_{kj}(\om)$. The unit transportation costs to $j$ are $c_{ij}$ for facility $i$ and $(q_{kj} - s_{kj} x''_k)$ for supplier $k$. Here, $q_{kj}$ represents the unit transportation cost from supplier $k$ to demand location $j$ without an outsourcing contract, i.e., $x''_k = 0$. This cost is reduced by $s_{kj} \geq 0$ if a contract is established, i.e., $x''_k = 1$. Let $h_i$ be the capacity of each facility $i \in I$. The first-stage feasible region is defined as $X = \{(x',x'') \in \{0,1\}^{n_1}\times\{0,1\}^{n_3}: \kappa_1^\top x' + \kappa_2^\top x'' \leq b\}$, where $\kappa_1 \in \R^{n_1}_+$ and $\kappa_2 \in \R^{n_3}_+$ are cost vectors associated with establishing facilities and outsourcing contracts, respectively. 

The formulation of risk-neutral BiFLP is given by
$
    \min_{x \in X} \sum_{\om \in \Om} p(\om) Q(x, \om)
$
where
\begin{subequations} \label{eq:flp-Q}
\begin{align}
    Q(x, \om) =
    \min_{y(\om) \in \R_+^{n_1 \times n_2}, \ u(\om) \in \R_+^{n_3 \times n_2}} \ & \sum_{i \in I, j \in J} c_{ij} y_{ij}(\om) + \sum_{k \in K, j \in J} (q_{kj} - s_{kj} x''_k) u_{kj}(\om) \label{eq:flp-obj2} \\
    \text{s.t.} \ 
    & \sum_{i \in I} y_{ij}(\om) + \sum_{k \in K} u_{kj}(\om) \geq d_{j}(\om), \quad \forall j \in J, \label{eq:flp-demand-con}\\
    & \sum_{j \in J} y_{ij}(\om) \leq h_{i} x'_i, \quad \forall i \in I, \label{eq:flp-capacity-con}
\end{align}
\end{subequations} 
for $\om \in \Om$. The objective~\cref{eq:flp-obj2} is to minimize the total cost of fulfilling demand, considering both transportation costs from facilities and outsourcing suppliers. Constraints~\cref{eq:flp-demand-con} ensure that demand at all locations $j \in J$ are satisfied. Constraint~\cref{eq:flp-capacity-con} for each $i \in I$ limits the total flow from facility $i$ by its capacity $h_i$. 

To generate test instances, we randomly place $(n_1 + n_2 + n_3)$ points on a $100 \times 100$ grid, representing potential facility locations, demand locations, and supplier locations. We consider four network sizes, with $(n_1, n_2, n_3)$ set to $(12, 40, 5),$ $(12, 40, 10),$ $(20, 60, 5),$ and $(20, 60, 10).$
The costs of building a facility $\kappa_{1i} = 5$ and contracting an outsourcing supplier $\kappa_{2k}=4$ for all $i \in I$ and $k \in K$, and the budget $b = (5n_1 + 4n_3)/4$. The cost of fulfilling demand at location $j \in J$ from supplier $k \in K$ consists of a fixed component and a distance-dependent component: specifically, $q_{kj} = \bar{c} + 2 v(k, j)$ and $s_{kj} = 
\bar{c} + 0.8 v(k,j)$, where $\bar{c}$ is a predetermined fixed cost, and $v(k, j)$ is the Euclidean distance between $k$ and $j$. Each facility at $i \in I$ has a capacity of $h_i = 100$. 
Demand data are generated using normal distributions. For instances with $n_3 = 5$, the mean demand $\bar{\mu}_j$ for each $j \in J$ is uniformly drawn from $\{40, 41, \dots, 80\}$, and for instances with $n_3 = 10$, it is drawn from $\{50, 51, \dots, 90\}$. In both cases, the standard deviation for each $j \in J$ is set to $\bar{\mu}_j / 4$.

For benchmark, we consider two approaches: an approach akin to the integer $L$-shaped method (IL) and a deterministic extensive formulation (DE). In IL, the recourse function is approximated using integer optimality cuts (see Proposition~2 in \citet{Laporte}). Unlike the standard integer $L$-shaped method, IL does not generate continuous optimality cuts, since the continuous relaxation of the recourse problem~\cref{eq:Q} does not provide dual information for generating such cuts.
For the min-min problem, DE is formulated as the following bilinear program: 
\begin{equation}
\begin{aligned}
    \min \bigg\{
    \sum_{\om \in \Om} p(\om) \Big(\sum_{i \in I, j \in J} c_{ij} y_{ij}(\om) + \sum_{k \in K, j \in J}(q_{kj} - s_{kj} x''_k) u_{kj}(\om)\Big) : \\ x \in X,\  \cref{eq:flp-demand-con}\text{--}\cref{eq:flp-capacity-con}, y(\om)\geq 0, u(\om) \geq 0, \om \in \Om\bigg\}.
\end{aligned}
\end{equation}
In our tests, this formulation is solved directly using Gurobi~9.5, with \textit{NonConvex} parameter set to 2.
In the standard $L^2$ method (denoted by $L^2$), we set $\sigma_i = 10^{5}$ for all $i \in I$. In the regularized $L^2$ method (denoted by $L^2$-R), we scale the objective by a factor of $10^{-2}$ to balance its magnitude with the regularization term. Additionally, we set $\gamma = 10^{-4}$ and $\bar{u}_i = 10^3$ for all $i \in I$ and define the regularization function as $R(\mu) = \sum_{i \in I} \mu_i$.

\begin{table}[t]
\centering
\begin{threeparttable}[t]
\fontsize{10}{12}\selectfont
\setlength{\tabcolsep}{4.9pt}
\caption{Results for BiFLP instances.}
\label{tab:min-min_result}
\begin{tabular}{llrrrrrrrr}
\toprule
\thead{2}{Instance} & \thead{2}{$L^2$} & \thead{2}{$L^2$-R} & \thead{2}{DE} & \thead{2}{IL} \\
\cmidrule(lr){1-2} 
\cmidrule(lr){3-4} \cmidrule(lr){5-6} \cmidrule(lr){7-8} \cmidrule(lr){9-10}
\thead{1}{$(n_1, n_2, n_3)$} & \thead{1}{$|\Om|$} & \thead{1}{Gap (\%)} & \thead{1}{Time (s)} & \thead{1}{Gap (\%)} & \thead{1}{Time (s)} & \thead{1}{Gap (\%)} & \thead{1}{Time (s)} & \thead{1}{Gap (\%)} & \thead{1}{Time (s)} \\
% \midrule
\specialrule{1pt}{2pt}{2pt}
$(12, 40, 5)$ &   10 & 0.0      & 2        & 0.0      & 2        & 0.0      & 2        & 0.0      & 20       \\
              &   50 & 0.0      & 6        & 0.0      & 7        & 0.0      & 27       & 0.0      & 78       \\
              &  100 & 0.0      & 12       & 0.0      & 12       & 0.0      & 95       & 0.0      & 142      \\
              &  200 & 0.0      & 23       & 0.0      & 25       & 0.0      & 344      & 0.0      & 296      \\
              &  500 & 0.0      & 49       & 0.0      & 54       & 0.0      & 1732     & 0.0      & 652      \\
              & 1000 & 0.0      & 131      & 0.0      & 131      & NA       & 3600+    & 0.0      & 1300     \\ \cmidrule(lr){1-10}
$(12, 40,10)$ &   10 & 0.0      & 18       & 0.0      & 33       & 0.0      & 20       & 0.0      & 469      \\
              &   50 & 0.0      & 60       & 0.0      & 67       & 0.0      & 320      & 0.0      & 1320     \\
              &  100 & 0.0      & 137      & 0.0      & 148      & 0.0      & 1168     & 0.0      & 2373     \\
              &  200 & 0.0      & 196      & 0.0      & 243      & NA       & 3600+    & 100.0    & 3600+    \\ \cmidrule(lr){1-10}
$(20, 60, 5)$ &   10 & 0.0      & 40       & 0.0      & 69       & 0.0      & 8        & 100.0    & 3600+    \\
              &   50 & 0.0      & 110      & 0.0      & 126      & 0.0      & 171      & 100.0    & 3600+    \\
              &  100 & 0.0      & 113      & 0.0      & 112      & 0.0      & 457      & 100.0    & 3600+    \\
              &  200 & 0.0      & 263      & 0.0      & 290      & 0.0      & 1399     & 100.0    & 3600+    \\ \cmidrule(lr){1-10}
$(20, 60,10)$ &   10 & 0.0      & 285      & 0.0      & 394      & 0.0      & 92       & 100.0    & 3600+    \\
              &   50 & 0.0      & 854      & 0.0      & 867      & 0.0      & 1697     & 100.0    & 3600+    \\
              &  100 & 0.0      & 1240     & 0.0      & 1422     & NA       & 3600+    & 100.0    & 3600+    \\
              &  200 & 5.0\tnote{*}      & 2423\tnote{*}     & 6.8\tnote{**}      & 2612\tnote{**}     & NA       & 3600+    & 100.0    & 3600+    \\
\bottomrule
\end{tabular}
\begin{tablenotes} \footnotesize
    \item [*] Average of the following results: (1) 0\% gap, 1909 s, (2) 0\% gap, 1756 s, and (3) 14.9\% gap, 3600 s.
    \item [**] Average of the following results: (1) 0\% gap, 2057 s, (2) 0\% gap, 2175 s, and (3) 20.5\% gap, 3600 s.
\end{tablenotes}
\end{threeparttable}
\end{table}

\Cref{tab:min-min_result} summarizes the test results for the BiFLP instances. Each row presents the average results of three instances with the same network structure, $(n_1, n_2, n_3)$, and the same number of scenarios, $|\Om|$. The columns labeled ``Gap (\%)'' and ``Time (s)'' report the optimality gap (in \%) and solution time (in seconds), respectively. The optimality gap results are marked as ``NA'' for instances where an algorithm failed to find both primal and dual bounds within the time limit.
The results show that $L^2$ outperforms the other approaches in terms of the computational efficiency. On average, $L^2$ is 17.6 times faster than IL and is 7.4 times faster than DE for the instances solved to optimality by all approaches. These factors increase to 21.9 and 8.8 times, respectively, when considering all instances. The IL showed poor scalability due to its limited capability in improving dual bounds; specifically, for the instances with $(n_1,n_2,n_3) = (20, 60, 5)$ and $(20, 60, 10)$, IL could not reduce optimality gaps within the time limit for all instances. We find that DE's performance is less sensitive to the network size than the others, but it is significantly affected by the number of scenarios. For the first instance category, with $(12, 40, 5)$ network and $10$ scenarios, DE and $L^2$ solved instances in similar solution times. However, as the number of scenarios increases to 500, DE's solution time increases by around 900 times, while $L^2$'s solution time increases only by around 26 times. When comparing the results from $L^2$ and $L^2$-R, the performance differences are minor in terms of solution time. The standard $L^2$ method is, on average, 1.1 times faster than the regularized $L^2$ method for instances where both methods solved to optimality.

\subsection{Bi-Parameterized (Min-Max) Network Interdiction Problem}\label{sec:compres-minmax}

Next, we consider a bi-parameterized (min-max) network interdiction problem (BiNIP). 
Consider a directed graph $\mathcal{G} = (\N, \A)$, where $\N = \{1, \dots, n_1\}$ is the set of nodes and $\A = \{1, \dots, n_2\}$ is the set of arcs in the graph. Resources that restrict each path in this network is indexed by $K = \{1, \dots, n_3\}$. Variable $x'_a \in \{0,1\}$ for each $a \in \A$ indicates whether arc $a \in \A$ is interdicted. For resources, $x''_k \in \{0,1\}$ for $k \in K$ represents whether interdiction occurs for resource $k$ or not. Let $x = (x', x'')$ represents a vector of all interdiction decision variables. The interdiction decisions are associated with costs $\kappa_a \in \R_+$ for arcs and $g_k \in \R_+$ for resources, and the total cost is constrained by budget $b \in \R_+$. The first-stage feasible region is defined as $ X = \{(x', x'') \in \{0,1\}^{n_2} \times \{0,1\}^{n_3} : \kappa^\top x' + g^\top x'' \leq b \}.$ 
Variable $y_a \in \{0,1\}$ for each $a \in \A$ represents whether the network user traverses arc $a$ ($y_a = 1$) or not ($y_a = 0$). Using random variable $\tilde{d}_a$, we represent the increase in arc length due to interdiction for each $a \in \A$. Its realizations are denoted by $d_{a}(\om)$ for $\om \in \Om$. The length of arc $a$ for scenario $\om$ becomes $(c_a + d_{a}(\om) x'_a)$, where $c_a$ is the nominal length of arc $a$ when not interdicted. The change in resource $k \in K$ after traversing arc $a \in \A$ is denoted by $r_{k a}$, and the threshold is denoted by $h_{k}$. When interdiction occurs, i.e., $x''_k = 1$, this threshold is adjusted by $s_{k}$. The formulation of BiNIP is given by
$
    \max_{x \in X} \sum_{\om \in \Om} p(\om) Q(x, \om)
$
where, for $\om \in \Om$,
\begin{subequations} \label{eq:nip-Q}
\begin{align}
    Q(x, \om) = \min_{y \in \{0,1\}^{n_2}} \ & \sum_{a \in \A} (c_a + d_{a}(\om) x'_a) y_a \label{eq:nip-obj} \\ 
    \text{s.t.} \ 
    & T y = q, \label{eq:nip-flow-balance} \\
    & \sum_{a \in \A} r_{k a} y_a \geq h_{k} + s_{k} x''_k, \quad \forall k \in K. \label{eq:nip-resource}
\end{align}
\end{subequations}
The first-stage problem aims to maximize the expected path length, with interdiction solutions restricted by the cardinality constraint in $X$. In the network user's problem, the objective function~\cref{eq:nip-obj} represents the length of the path. Constraints~\cref{eq:nip-flow-balance} enforce the balance of incoming and outgoing flows for each node; $T \in \{-1,0,1\}^{n_1 \times n_2}$ is the node-arc incidence matrix, where $T_{i a} = 1$ if arc $a \in \A$ leaves node $i \in \N$, $T_{i a} = -1$ if arc $a$ enters node $i$, and $T_{i a} = 0$ otherwise. Also, $q \in \{-1,0,1\}^{n_2}$ is a vector where $q_i = 1$ if $i \in \N$ is the source node, $q_i = -1$ if $i$ is the sink node, and $q_i = 0$ otherwise. Constraints~\cref{eq:nip-resource} are the resource constraints. It is important to note that the integral restrictions on $y$ are necessary---unlike in the conventional shortest path problem---since the resource constraints may eliminate the integral property of the feasible region.

In our experiments, we use randomly generated instances based on the instances from \citet{NguyenSmith}. We utilize their data on network topology, arc lengths, and deterministic penalty lengths. We consider two categories of their instances: $20$-node and $40$-node instances, with 10 instances in each category. The number of arcs varies in $[59, 78]$ for the $20$-node instances, and $[277,306]$ for the $40$-node instances. We extend these instances by introducing random penalty lengths, which are sampled from a uniform distribution over the interval $[\bar{d}_a - o_a, \bar{d}_a + o_a]$ for each arc $a \in \A$, where $\bar{d}_a$ is the deterministic penalty length from \citet{NguyenSmith}. The offset $o_a = 3$ for the $20$-node instances and $o_a = 4$ for the $40$-node instances.
For the $20$-node instances, we set budget $b=4$, arc interdiction costs $\kappa_a = 1$ for all $a \in \A$, and resource interdiction costs $g_k = 2$ for all $k \in K$. We consider three resources ($K = \{1, 2, 3\}$). The resource consumption parameter $r_{k a}$ is randomly drawn from $\{1, 2\}$ for each $k \in K$ and $a \in \A$. The threshold vector $h = (5, 4, 3)^\top$, and the penalty vector $s = (1, 2, 3)^\top$. Similarly, for the $40$-node instances, budget $b$ is set to $5$ with the same arc and resource interdiction costs. We consider four resources ($K = \{1, 2, 3, 4\}$), where $r_{k a}$ is randomly drawn from $\{1, 2, 3\}$ for each $k \in K$ and $a \in \A$. We set the threshold vector $h = (10, 8, 6, 4)^\top$ and the penalty vector $s = (1, 3, 5, 7)^\top$.

To benchmark the proposed approach, we consider IL for BiNIP. Note that DE is not applicable to BiNIP due to its min-max form. In preliminary tests, we found minor differences between the outcomes of $L^2$ and $L^2$-R. Therefore, we report only the results obtained by $L^2$ for BiNIP in \Cref{tab:min-max_result}. For all tests, the parameters $\sigma_k$ are set to $10^{3}$ for $k \in K$.
\begin{table}[t]
\centering
% \fontsize{10}{12}\selectfont
% \setlength{\tabcolsep}{4.8pt}
\caption{Results for BiNIP instances.}
\label{tab:min-max_result}
\begin{tabular}{llllrrrrr}
\toprule
\thead{4}{Instance} & \thead{2}{$L^2$} & \thead{3}{IL}  \\ \cmidrule(lr){1-4} \cmidrule(lr){5-6} \cmidrule(lr){7-9}
$n_1$ &       $n_2$ & $n_3$ & $|\Om|$ & Gap (\%)         & Time (s)        & Gap (\%)    & RelGap (\%) & Time (s)  \\ \specialrule{1pt}{1.5pt}{1.5pt}
   20 &   $[59,78]$ &     3 & 10      & 0.0        & 0.6         & 100.0 & 3.6   & 3600+ \\ 
      &             &       & 20      & 0.0        & 1.4         & 100.0 & 4.0   & 3600+ \\
      &             &       & 50      & 0.0        & 3.6         & 100.0 & 6.8   & 3600+ \\
      &             &       & 100     & 0.0        & 7.2         & 100.0 & 6.4   & 3600+ \\
      &             &       & 500     & 0.0        & 41.7        & 100.0 & 11.2  & 3600+ \\
      &             &       & 1000    & 0.0        & 60.9        & 100.0 & 9.5   & 3600+ \\ \cmidrule(lr){1-9}
   40 & $[277,306]$ &     4 & 10      & 0.0        & 5.8         & 100.0 & 3.6   & 3600+ \\
      &             &       & 20      & 0.0        & 10.8        & 100.0 & 4.0   & 3600+ \\
      &             &       & 50      & 0.0        & 31.6        & 100.0 & 4.3   & 3600+ \\
      &             &       & 100     & 0.0        & 63.6        & 100.0 & 5.3   & 3600+ \\ 
\bottomrule
\end{tabular}
\end{table}
The test results are summarized in \Cref{tab:min-max_result} where each row presents the average results for 10 instances. For each column labeled ``$L^2$'' or ``IL'', we report the optimality gap (in \%) under ``Gap (\%)'', and the solution time (in seconds) under ``Time (s)''. The results under the ``RelGap (\%)'' column represent the relative gaps between the primal bounds obtained by IL to the optimal objective values.
The results show that the $L^2$ method outperforms IL across all instances. The IL was unable to reduce the dual bounds for all 100 instances within the time limit of 3600 seconds, resulting in 100\% optimality gaps, while $L^2$ found optimal solutions for all instances within 23 seconds on average. When comparing primal bounds, IL produced primal bounds that were, on average, $5.7\%$ worse than those obtained by the $L^2$ method, even if the former spent 158 times more computational time.

\subsection{Distributionally Robust Facility Location with Decision-Dependent Ambiguity Set} \label{sec:test_result-dddro}

We now consider a distributionally robust two-stage facility location problem under decision-dependent demand uncertainty (DRFLP), which is modified from the problem presented in \citet{Yu}. In DRFLP, we assume that facility locations impact service accessibility, thereby affecting the realizations of random demand. While the exact distribution of demand is unknown, we model its moments as functions of chosen locations to express how demand depends on these decisions. Specifically, these functions in our model capture the dynamics that locating a facility closer to a demand point results in a higher mean demand, compared to placing it farther away. To determine robust decisions under this demand uncertainty and distributional ambiguity, we employ a DRO model, where the ambiguity set is defined using these moment functions.

Let $I:=\{1, \dots, n_1\}$ be the set of potential facility locations, and $J := \{1, \dots, n_2\}$ the set of demand locations. The binary variable $x_i$ indicates whether a facility is established at location $i \in I$, with the total number of establishments constrained by budget $b$. The random demand at each $j \in J$ is denoted by $\tilde{d}_j$, and its realization is denoted by $d_{j}(\om)$ for $\om \in \Om$. The variables $y_{ij}$ represents the flow from facility $i \in I$ to demand location $j \in J$, associated with unit flow cost $c_{ij}$. If the demand at location $j \in J$ is not satisfied, each unit of unsatisfied demand incurs a penalty cost $q_j$. Each facility $i$ has capacity $h_i$, restricting the total outgoing flow. DRFLP can then be formulated as
\begin{equation}
    \min_{x \in X} \max_{p \in \mathcal{P}(x)} \sum_{\om \in \Om} p(\om) Q(x, \om)
\end{equation}
where $X =  \{x \in \{0,1\}^{n_1} : \sum_{i \in I} x_i \leq b\}$ and for $\om \in \Om$,
\begin{subequations} \label{eq:drflp-Q}
\begin{align}
    Q(x, \om) = \min_{y \in \R_+^{n_1 \times n_2}, \ u \in \R_+^{n_2}} \ &
        \sum_{i \in I, j \in J} {c_{ij} y_{ij}} + \sum_{j \in J} q_{j} u_{j}  \label{eq:drflp-Q-obj} \\
        \text{s.t.} \ 
        & \sum_{i \in I} y_{ij} + u_j \geq d_{j}(\om), \quad \forall j \in J \label{eq:drflp-Q-b} \\
        & \sum_{j \in J} y_{ij} \leq h_i x_i, \quad \forall i \in I. \label{eq:drflp-Q-c}
\end{align}
\end{subequations} 
The objective of the recourse problem is to minimize the sum of transportation and penalty costs as described in \cref{eq:drflp-Q-obj}. Constraints~\cref{eq:drflp-Q-b} ensure that the total flow into demand location $j \in J$, along with the amount $u_j$, satisfies the demand $d_{j}(\om)$. Constraints~\cref{eq:drflp-Q-c} limits the total flow from each facility $i \in I$ by its capacity $h_i$ if the facility is established (i.e., $x_i=1$), or to zero otherwise.
The moment-matching ambiguity set $\mathcal{P}(x)$ is defined as follows.
\begin{equation}
\begin{aligned}
    \mathcal{P}(x) = \bigg\{ & p \in \R_+^{|\Om|} :  
        \sum_{\om \in \Om} p(\om) d_{j}(\om) \in \Big[M_j(x) - \epsilon^M_j, M_j(x)+\epsilon^M_j\Big], \  \forall j \in J, \\
        & \sum_{\om \in \Om} p(\om) (d_{j}(\om))^2 \in \Big[S_j(x) \underline{\epsilon}^S_j, S_j(x) \bar{\epsilon}^S_j \Big], \  \forall j \in J, 
        \ \sum_{\om \in \Om} p(\om) = 1
        \bigg\}.
\end{aligned}
\end{equation}
This ambiguity set consists of bounding constraints on the first and second moments of $\tilde{d}_j$ for $j \in J$ under the probability distribution $p$. The parameters are defined as follows:
\begin{equation*}
    M_j(x) = \bar{\mu}_j \Big(1 + \sum_{i \in I} \rho^{M}_{ij} x_i \Big), \ \text{ and }\ 
    S_j(x) = \big(\bar{\mu}_j^2 + \bar{\sigma}^2_j \big) \Big(1 + \sum_{i \in I} \rho^{S}_{ij} x_i \Big), \quad \forall j \in J.
\end{equation*}
Here, $\bar{\mu}_j$ and $\bar{\sigma}_j$ denote the baseline first and second moments for $j \in J$. The parameters $\rho^M_{ij} = e^{-v(i,j) / 25}$ and $\rho^S_{ij} = e^{-v(i,j) / 50}$ represent the effect of establishing facility at $i$ on the moments of $\tilde{d}_j$, where $v(i,j)$ is the Euclidean distance between locations $i \in I$ and $j \in J$.  
We set $\epsilon^M_j = 25$, $\underline{\epsilon}^S_j = 0.1$, and $\bar{\epsilon}^S_j = 1.9$ for all $j \in J$.

To generate test instances, we place $(n_1 + n_2)$ random locations on a $100\times100$ grid for facility and demand locations. The transportation cost $c_{ij}$ is set to the Euclidean distance $v(i, j)$ from location $i \in I$ to location $j \in J$. We set the capacity $h_i = h = 500$ for $i \in I$. The mean values $\bar{\mu}_j$ for $j \in J$ are uniformly sampled from $\{\underline{d}_j, \underline{d}_j + 1, \dots, \bar{d}_j\},$ where $\underline{d}_j$ is the nearest integer to $(0.7 \times b \times h / n_2)$ and $\bar{d}_j$ is the nearest integer to $\round(b \times h / n_2)$. Then, we sample the demand realizations $d_{j}(\om)$ for $\om \in \Om$ from $N(\bar{\mu}_j, 0.8 \bar{\mu}_j)$, truncated within $[1, 300]$, where $0.8 \bar{\mu}_j$ is the standard deviation.  

In the $L^2$ methods for DRFLP, we include an additional step that determines the worst-case distribution within the ambiguity set after solving subproblems. In particular, this step involves solving the \textit{distribution separation problem}, defined as follows:
\begin{equation}
    \max_{p \in \mathcal{P}(z), z \in X} \sum_{\om \in \Om} p(\om) \hat{Q}^l(x^{l}, \om) + \sum_{i \in I} \lambda^l_i z_i
\end{equation}
where $\lambda^l_i = \sigma_i (2x^{l}_i - 1)$ in \Cref{alg:min-max-L2} or $\lambda^l_i = \mu^{l}_i (2x^{l}_i - 1)$ in \Cref{alg:min-max_alt}, respectively.
We denote the worst-case distribution identified by solving the distribution separation problem in iteration $l$ by $p^l = (p^l(\om))_{\om \in \Om}^\top$. Using $p^l$, we evaluate the objective $\phi^l = f(x^{l}) - \sigma^\top x^{l} + \sum_{\om \in \Om} p^l(\om) \hat{Q}^l(x^{l}, \om)$ and obtain an optimality cut in the following form:
\begin{equation*}
    \theta \geq \sum_{\om \in \Om} p^l(\om) \Big( 
        \sum_{j \in J} d_{j}(\om) \pi^l_{j}(\om) - \sum_{i \in I} h_i x_i \nu^l_{i}(\om) \Big)
        + \sum_{i \in I} \sigma_i  (2x_i - 1) z^l_{i},
\end{equation*}
where $(\pi^l(\om), \nu^l(\om))$ is an optimal dual solution for the subproblem, defined by the LP dual of \cref{eq:drflp-Q}, at scenario $\om$. 
To reduce the computational burden, we fix $z$ to $x^{l}$ in the distribution separation problem. 
If the ambiguity set is empty for a given current solution, i.e., $\P(x^l) = \emptyset$, then we add the following feasibility cut: 
$
    \sum_{i \in I | x^{l}_i = 0} x_i + \sum_{i \in I | x^{l}_i = 1} (1- x_i) \geq 1.
$

\begin{table}[t]
\centering
% \fontsize{10}{12}\selectfont
% \setlength{\tabcolsep}{4.5pt}
\caption{Results from $L^2$-R and DA for the DRFLP instances.
}\label{tab:drflp-result}
\begin{tabular}{llrrrrrr}
\toprule
\thead{2}{Instance} & \thead{3}{$L^2$-R} & \thead{3}{DA} \\ \cmidrule(lr){1-2} \cmidrule(lr){3-5} \cmidrule(lr){6-8}
\thead{1}{$(n_1, n_2, b)$} & \thead{1}{$|\Om|$} & \thead{1}{Gap (\%)} & \thead{1}{Time (s)} & \thead{1}{\#OptCut} & \thead{1}{Gap (\%)} & \thead{1}{Time (s)} & \thead{1}{\#OptCut} \\
\specialrule{1pt}{1.5pt}{1.5pt}
(15, 20, 4)    & 500     & 0.0      & 34       & 101      & 0.0      & 39       & 47013    \\
      & 1000    & 0.0      & 37       & 86       & \multicolumn{3}{c}{Unbounded}  \\
      & 2000    & 0.0      & 48       & 61       & 0.0      & 76       & 94087    \\
      & 5000    & 0.0      & 225      & 118      & 0.0      & 621      & 583829   \\
      & 10000   & 0.0      & 578      & 140      & 0.0      & 3589     & 1552667  \\ \cmidrule(lr){1-8}
(20, 20, 4)    & 500     & 0.0      & 27       & 104      & 0.0      & 38       & 48476    \\
      & 1000    & 0.0      & 65       & 143      & 0.0      & 212      & 131379   \\
      & 2000    & 0.0      & 144      & 177      & 0.0      & 565      & 322976   \\
      & 5000    & 0.0      & 515      & 222      & 0.0      & 2178     & 849847   \\
      & 10000   & 0.0      & 844      & 156      & 22.6     & 3600+    & 1298870  \\ \cmidrule(lr){1-8}
(30, 20, 4)    & 500     & 0.0      & 379      & 937      & \multicolumn{3}{c}{Unbounded}  \\
      & 1000    & 0.0      & 536      & 885      & 21.4     & 3600+    & 644380   \\
      & 2000    & 0.0      & 320      & 276      & 0.0      & 3273     & 514143   \\
      & 5000    & 0.0      & 1228     & 457      & \multicolumn{3}{c}{Unbounded}  \\
      & 10000   & 0.0      & 1683     & 286      & 100.0    & 3600+     & 932400   \\ \cmidrule(lr){1-8}
(40, 20, 4)    & 500     & 0.0      & 492      & 989      & \multicolumn{3}{c}{Unbounded}  \\
      & 1000    & 0.0      & 860      & 1132     & \multicolumn{3}{c}{Unbounded}  \\
      & 2000    & 0.0      & 677      & 459      & \multicolumn{3}{c}{Unbounded}  \\
      & 5000    & 7.9      & 3600+    & 1125     & 72.8     & 3600+    & 805440   \\
      & 10000   & 0.0      & 1926     & 259      & \multicolumn{3}{c}{Unbounded}  \\ 
\bottomrule
\end{tabular}
\end{table}

\Cref{tab:drflp-result} presents the test results comparing the performance of the $L^2$ method and the dual-based approach (DA). Here, we focus on the regularized $L^2$ method ($L^2$-R, \Cref{alg:min-max_alt}), as it showed a better performance in our preliminary tests for DRFLP. For $L^2$-R, we set $\gamma = 10^{-1}$ and $\bar{u}_i = 10^{4}$ for all $i \in I$ and define the regularization function as $R(\mu) = \norm{\mu}^2_2$. The DA solves the dual reformulation~\cref{eq:dddro-dual-reform} of DRFLP using the decomposition algorithm presented in \cite{Yu}. We consider different instance settings by varying the number of facility locations $n_1$ in $\{15, 20, 30, 40\}$, and the number of scenarios $|\Om|$ in $\{500, 1000, 2000, 5000, 10000\}$. We fix the number of demand locations $n_2 = 20$ and the budget $b = 4$. 

Each row in \Cref{tab:drflp-result} reports the optimality gap ``Gap (\%)'', solution time in seconds ``Time (s)'', and the number of optimality cuts ``\#OptCut'' for each instance. Out of the total 20 instances, $L^2$-R solves 19 instances to optimality within the time limit, while DA solves only nine instances. DA reports ``unbounded" for instances where the ambiguity sets are not relative complete, whereas $L^2$-R successfully handles these instances. When comparing solution time, $L^2$-R is on average 5.3 times faster than DA for instances that are solved to optimality by both methods. As the number of scenarios increases, the difference in the number of optimality cuts increases, as DA requires significantly more cuts. On average, the $L^2$ method achieves optimality by adding only $0.04\%$ of the number of cuts generated/used by DA.

\begin{table}[t]
\centering
% \fontsize{10}{12}\selectfont
\caption{Results from $L^2$ and $L^2$-R for DRFLP instances.} \label{tab:drflp-result-l2-method}
\begin{tabular}{llrrrr}
\toprule
\thead{2}{Instance} & \thead{2}{$L^2$} & \thead{2}{$L^2$-R} \\ \cmidrule(lr){1-2} \cmidrule(lr){3-4} \cmidrule(lr){5-6}
\thead{1}{$(n_1, n_2, b)$} & \thead{1}{$|\Om|$} & \thead{1}{Gap (\%)} & \thead{1}{Time (s)} & \thead{1}{Gap (\%)} & \thead{1}{Time (s)} \\ \specialrule{1pt}{1.5pt}{1.5pt}
(20, 20, 4) & 100   & 0.0  & 162  & 0.0 & 14  \\
            & 300   & 0.0  & 324  & 0.0 & 36  \\
            & 500   & 0.0  & 620  & 0.0 & 27  \\
            & 1000  & 0.0  & 1125 & 0.0 & 65  \\
            & 2000  & 0.0  & 2405 & 0.0 & 144 \\
            & 5000  & 51.7 & 3600+ & 0.0 & 515 \\
            & 10000 & 63.1 & 3600+ & 0.0 & 844 \\ \bottomrule
\end{tabular}
\end{table}

Next, we compare the performance of the standard $L^2$ method ($L^2$, \Cref{alg:min-max-L2}) with the regularized $L^2$ method ($L^2$-R, \Cref{alg:min-max_alt}) in \Cref{tab:drflp-result-l2-method}. Parameters $\sigma_i$ for $i \in I$ are set to $10^3$. Test instances have $(n_1, n_2, b) = (20, 20, 4)$ and the number of scenarios $|\Om| \in \{100, 300, 500, 1000, 2000, 5000, 10000\}$. The results show that $L^2$-R is computationally efficient than $L^2$ for solving the DRFLP instances. The $L^2$-R solves all instances, while $L^2$ is unable to solve the instances with $5000$ and $10000$ scenarios. For instances where both methods solve to optimality, $L^2$-R is, on average, 16.2 times faster than $L^2$ method.

\section{Conclusion} \label{sec:conclusion}

In this paper, we introduced Lagrangian-integrated $L$-shaped ($L^2$) methods for solving bi-parameterized two-stage stochastic (min-max and min-min) integer programs (BTSPs), which are applicable to interdiction models and a wide range of optimization problems with decision-dependent uncertainty. 
For cases where the first-stage decisions are pure binary, we developed two exact algorithms applied to the min-min and min-max BTSPs. Additionally, we proposed a regularization-augmented method to address BTSPs with mixed-integer first-stage decisions. 
To evaluate the $L^2$ method's efficiency, we conducted extensive numerical tests in various settings: 
\textit{(a)} a bi-parameterized (min-min) facility location problem; \textit{(b)} a bi-parameterized (min-max) network interdiction problem; and \textit{(c)} a distributionally robust facility location problem with a decision-dependent ambiguity set, where nonemptiness is not guaranteed for all first-stage decisions---that is, \textit{non-relatively complete} ambiguity set. The results demonstrated the superior efficiency of the $L^2$ method compared to benchmark approaches. Specifically, our approach converged to optimal solutions for all tested instances of the bi-parameterized network interdiction problem within 23 seconds on average, while the benchmark method failed to converge for any instance within 3600 seconds. The $L^2$ method achieved optimal solutions, on average, 18.4 times faster for the bi-parameterized facility location problem (with known distribution). For the decision-dependent DRO setting, the $L^2$ method effectively solved instances with the non-relatively complete ambiguity sets and achieved solutions 5.3 times faster than an existing dual-based approach. 

\subsection*{Acknowledgements} This research is partially funded by National Science Foundation Grant CMMI-1824897 and 2034503, and Commonwealth Cyber Initiative grants which are gratefully acknowledged.

\newpage
\begin{appendices}

\section{BTSP-Based Reformulations for Stochastic Optimization and Interdiction Problems} \label{apx:btsp-reform}

In this section, we present the details of BTSP-based reformulations for decision-dependent (risk-averse) stochastic optimization and generalized interdiction problems, introduced in \Cref{sec:introduction}. For clarity, we restate the single-parameterized recourse function $Q^s$, which will be used throughout this section, as follows:
$$Q^s(x, \om) = \min \{ q(\om)^\top y: W(\om) y = r(\om) - T(\om) x,\ y \in \Y \}.$$

\subsection{Risk-Averse Stochastic Optimization with CVaR-Based Decision-Dependent Uncertainty}

We consider the following formulation of two-stage stochastic programs with decision-dependent probabilities and CVaR measure:
\begin{equation} \label{eq:cvar}
\begin{aligned}
    \min_{x \in X} 
    \Big\{
    f(x) + \cvar_\alpha \big(Q^s(x, \om), p(x) \big)
    \Big\}.
\end{aligned}
\end{equation}
Using the linear programming formulation of CVaR~\cref{eq:def-cvar}, we rewrite this problem as follows:
\begin{subequations} \label{eq:cvar-reform0}
\begin{align} 
    \min \ & f(x) + \eta + \frac{1}{1-\alpha} \sum_{\om \in \Om} p(x, \om) \nu(\om) \\
    \text{s.t.}\ 
    & \nu(\om) \geq Q^s(x, \om) - \eta, \quad \forall \om \in \Om \\
    & x \in X, \ \eta \in \R, \ \nu(\om) \in \R_+, \ \forall \om \in \Om.
\end{align}
\end{subequations}
Suppose that $p(\cdot, \om)$ for each scenario $\om \in \Om$ is an affine function.
The above formulation can be addressed using a dual decomposition-based approach, as described in \cite{Schultz}. 
However, this approach only yields a lower bound on the optimal objective value, thereby leading to a duality gap potentially. Additionally, as discussed in \Cref{sec:introduction}, applying the dual decomposition method to this formulation results in mixed-integer nonconvex subproblems, which can impose a significant computational burden.

Now, we present a reformulation of \cref{eq:cvar} into the form of the min-max BTSP \eqref{eq:btsp}. By taking the dual of \cref{eq:def-cvar}, we can obtain the dual representation of $\cvar_\alpha (R (x, \om), p(x))$:
\begin{equation} \label{eq:cvar-dual}
    \max \bigg\{
        \frac{1}{1-\alpha}\sum_{\om \in \Om} \gamma(\om) Q^s(x, \om) : \sum_{\om \in \Om} \gamma(\om) = 1-\alpha, \ 0 \leq \gamma(\om) \leq p(x, \om), \ 
        \forall \om \in \Om
    \bigg\}.
\end{equation}
Let $Q^s(x, \om) = \max\{ \pi(\om)^\top (r(\om) - T(\om)x) : W(\om)^\top \pi(\om) \leq q(\om), \ \pi(\om) \in \R^{m_2}_+ \}$ be the dual formulation of the recourse problem with a convexified feasible region for each scenario $\om$. By incorporating this into the dual representation of $\cvar$, we have
\begin{subequations} \label{eq:cvar-dual-reform1}
\begin{align} 
    \cvar_\alpha (Q^s(x, \om), p(x)) = 
    \max \ & \frac{1}{1-\alpha}\sum_{\om \in \Om} \gamma(\om) \pi(\om)^\top(r(\om) - T(\om) x) \\
    \text{s.t.} \ 
    & \sum_{\om \in \Om} \gamma(\om) = 1-\alpha \\
    & \ \gamma(\om) \leq p(x, \om), \quad \forall \om \in \Om\\
    & W(\om)^\top \pi(\om) \leq q(\om), \quad \forall \om \in \Om \label{eq:cvar-dual-reform1-d}\\
    & \pi(\om) \in \R^{m_2}_+, \ \gamma(\om) \in \R_+, \quad \forall \om \in \Om.
\end{align}
\end{subequations}
For any $\gamma(\om) \geq 0$, the inequality $(W(\om)^\top \pi(\om) \leq q(\om))$ holds if and only if $(W(\om)^\top \gamma(\om) \pi(\om) \leq q(\om) \gamma(\om))$. Therefore, we can replace \cref{eq:cvar-dual-reform1-d} with $(W(\om)^\top \gamma(\om) \pi(\om) \leq q(\om) \gamma(\om))$ for each $\om \in \Om$. Next, we introduce a decision vector $\tau(\om) \in \R^{m_2}_+$ to substitute for $\gamma(\om) \pi(\om)$, as for any $\tau(\om) \geq 0$, there exist $\gamma(\om) \geq 0$ and $\pi(\om)$ such that $\tau(\om) = \gamma(\om) \pi(\om)$, and vice versa. Thus, we can reformulate \cref{eq:cvar-dual-reform1} as follows:
\begin{subequations} \label{eq:cvar-dual-reform2}
\begin{align} 
    \cvar_\alpha (Q^s(x, \om), p(x)) = 
    \max \ & \frac{1}{1-\alpha}\sum_{\om \in \Om} \tau(\om)^\top(r(\om) - T(\om) x) \\
    \text{s.t.} \ 
    & \sum_{\om \in \Om} \gamma(\om) = 1-\alpha \\
    & \ \gamma(\om) \leq p(x, \om), \quad \forall \om \in \Om\\
    & W(\om)^\top \tau(\om) - q(\om) \gamma(\om) \leq 0, \quad \forall \om \in \Om \\
    & \tau(\om) \in \R^{m_2}_+, \ \gamma(\om) \in \R_+, \quad \forall \om \in \Om,
\end{align}
\end{subequations}
which is in the form of the min-max problem~\cref{eq:btsp}.

\subsection{Two-Stage DRO with Decision-Dependent Ambiguity Set}

A two-stage decision-dependent DRO problem is defined as follows:
\begin{equation} \label{eq:dddro}
    \min_{x \in X} 
    \Big\{
        f(x) + \max_{p \in \P(x)} \E_p \big[Q^s(x, \om) \big]
    \Big\},
\end{equation}
where $\P(x)$ is an ambiguity set that depends on $x$, and $Q^s(x, \om) = \min\{q(\om)^\top y : W(\om) y = r(\om) - T(\om) x, \ y \in \Y\}$ for $\om \in \Om$. An example of an ambiguity set is \textit{moment-matching ambiguity set}, where the distribution’s moments match the known moment information. Specifically, let $\zeta(\om) := (\zeta_1(\om), \zeta_2(\om), \dots, \zeta_t(\om))^\top \in \R^t$ be moment functions on $\Om$. Then, a decision-dependent moment-matching ambiguity set is defined as:
\begin{equation} \label{eq:moment-ambiguity}
    \P(x) := \bigg\{
        p \in \R_+^{|\Om|} : 
        l(x) \leq \sum_{\om \in \Om} p(\om) \zeta(\om) \leq u(x), \underline{p}(x) \leq p \leq \bar{p}(x), \ \sum_{\om\in\Om} p(\om) = 1
    \bigg\}.
\end{equation}
Here, $\underline{p}(x):\R^{n_1} \to \R^{|\Om|}$, $\bar{p}(x):\R^{n_1} \to \R^{|\Om|}$, $l(x):\R^{n_1} \to \R^{t}$, and $u(x):\R^{n_1} \to \R^{t}$ are predetermined functions that specify lower and upper bounds for a given $x$. The DRO problem~\cref{eq:dddro} can be seen as a min-max formulation, where the recourse is associated with the decision $p$. We call the problem~\cref{eq:dddro} has \textit{relatively complete} ambiguity set, if $\P(x) \neq \emptyset$ for all $x \in X$.
It is important to note that the $L^2$ methods can address the DRO problem~\cref{eq:dddro} even in the absence of the relatively complete ambiguity set. Specifically, in the $L^2$ methods, by utilizing certain types of cuts, we can cut off infeasible solutions $x$ where $\P(x) = \emptyset$ while running the algorithm. For instance, if the ambiguity set is empty for a solution $\hat{x}$, then we can cut it off from the feasible region $X$ by adding the following cut:
\begin{equation}
    \sum_{i \in I | \hat{x}_i = 0} x_i + \sum_{i \in I | \hat{x}_i = 1} (1- x_i) \geq 1.
\end{equation}
This presents a distinct advantage of our approach when compared to an existing approach in the literature that is based on duality results. 

In the dual-based approach for \cref{eq:dddro}, we dualize the inner maximization problem, using strong duality of some special types of ambiguity sets, to derive a single-level reformulation, so-called a dual reformulation (e.g., see \citep{Basciftci, Luo2020, Yu}).
For the moment ambiguity set~\cref{eq:moment-ambiguity}, the dual reformulation of the DRO model~\cref{eq:dddro} is given by
\begin{subequations} \label{eq:dddro-dual-reform}
\begin{align}
    \min \ &
        f(x) - \underline{\alpha}^\top l(x) + \bar{\alpha}^\top u(x) - \underline{\beta}^\top \underline{p}(x) + \bar{\beta}^\top \bar{p}(x) \\
    \text{s.t.} \ 
    & (-\underline{\alpha} + \bar{\alpha})^\top \zeta(\om) - \underline{\beta}(\om) + \bar{\beta}(\om) \geq Q(x, \om), \quad \forall \om \in \Om \label{eq:dddro-dual-reform-con} \\
    & (\underline{\alpha}, \bar{\alpha}, \underline{\beta}, \bar{\beta}) \in \R^{t}_+ \times \R^{t}_+ \times \R^{|\Om|}_+ \times \R^{|\Om|}_+,\ x \in X,
\end{align}
\end{subequations}
where $\underline{\alpha}, \bar{\alpha}, \underline{\beta} = (\underline{\beta}(\om))_{\om \in \Om}^\top,$ and $\bar{\beta}= (\bar{\beta}(\om))_{\om \in \Om}^\top$ are dual multipliers for the constraints in \cref{eq:moment-ambiguity}. 
However, this dual approach presents computational challenges in practice. First, the dual reformulations rely on relatively complete ambiguity sets, which may be impractical, as demonstrated by our computational results in \Cref{sec:test_result-dddro}. One might consider using the penalty method---introducing penalties to address violated solutions---but it fails in the DRO problem~\cref{eq:dddro} due to its min-max structure. In the inner maximization problem, penalties must be applied negatively; however, these negative penalties can promote violations in the outer minimization problem, rather than prevent them. Another challenge is the scalability of the problem. Specifically, $Q(x, \om)$ in \cref{eq:dddro-dual-reform-con} is typically approximated using valid cuts. As the number of scenarios increases, the decomposed problems become increasingly difficult to solve due to the growing number of cuts. The nonconvex terms in the objective function also present further challenges in solving the decomposed problems. We note that the scalability issue is not limited to this specific type of ambiguity sets. When considering an ambiguity set defined using Wasserstein metric, so-called \textit{Wasserstein ambiguity set}, it is required to add $|\Om|\times|\Om|$ cuts in each iteration, readily resulting in a substantially large subproblem; e.g., see the algorithm presented in \citet{duque2020}.

\subsection{Bi-Parameterized Stochastic Network Interdiction Problem}

Recall that the generic formulation of stochastic interdiction problems is given by
\begin{equation*}
    \min_{x \in X} \E\Big[ \max_{y \in Y(x, \om)} f(x, y, \om) \Big].
\end{equation*}
Most studies investigating these problems typically assume that the interdiction $x$ affects either the objective function $f(x, y, \om)$ or the network user's feasible set $Y(x, \om)$, but not both, to derive efficient solution approaches (see, e.g.,  \citet{cormican1998,Morton2007,NguyenSmith}).

The min-max form of BTSP~\cref{eq:btsp} generalizes stochastic network interdiction problems by relaxing this assumption, thereby allowing for the modeling of more realistic situations. For instance, consider a network user seeking the shortest path on a directed graph $\G=(\N, \A)$, where feasible paths are subject to resource constraints. These resources can represent any values that may change during travel along an arc, such as travel time, fuel consumption, or load weight. These constraints ensure that the total resource usage along a path either meets or falls within specified thresholds. In this context, the interdictor may disrupt the network user's overall resource system, making it more challenging to satisfy the resource constraints. The resource constraint can be expressed by the following form, where the right-hand side depends on the interdiction decision $x''_k$ for each resource $k \in K$:
\begin{equation}
    \sum_{a \in \A} r_{\om k a} y_{\om a} \leq h_{\om k} + s_{\om k} x''_k,
\end{equation}
where, for scenario $\om$, $r_{\om k a} \in \R$ represents the resource change on arc $a$, $h_{\om k} \in \R$ expresses the nominal resource threshold, and $s_{\om k} \in \R$ denotes the impact of interdiction on the threshold.

Note that the interpretation of these constraints is not limited to resource contexts; for instance, in a surveillance coverage scenario (for the network user), the interdictor could force the network user to pass through specific nodes or arcs, causing a detour to the surveillance destination.
Incorporating these constraints into the network user's problem introduces integral restrictions on variables. While it is well-known that the feasible region of the conventional shortest path problem is integral---allowing the problem to be solved using its continuous relaxation without compromising optimality---introducing resource constraints eliminates this integral property. As a result, the network user's problem is required to have the integral restrictions on decision variables, i.e., $y_{\om a} \in \{0,1\}$ for all $a \in \A$ and $\om \in \Om$. Refer to \Cref{sec:compres-minmax} for results of our computational experiments for bi-parameterized stochastic network interdiction problem.

\section{Proofs} \label{apx:proofs}

\subsection{Proof of \texorpdfstring{\Cref{thm:strong_duality_Q_om}}{Theorem~\ref{thm:strong_duality_Q_om}}}
\begin{proof}{Proof}
    Clearly, $Q(x, \om) \leq D(x, \om)$ for any $\om\in\Om$ and $x \in X$. Therefore, to prove the statement, it suffices to show that $Q(x, \om) \geq D(x, \om)$, for all $\om \in \Om$ and $x \in X$.
    Fix $\om$, and
    let $Z(\om)$ denote the feasible set of the Lagrangian relaxation~\cref{eq:lag_Q_om}. Also, we let $\tilde{T}(\om) \in \Qu^{\tilde{m} \times n_1}, \tilde{W}(\om) \in \Qu^{\tilde{m} \times n_2},$ and $\tilde{r}(\om) \in \Qu^{\tilde{m}}$ such that $\conv(Z(\om)) = \{(y, z) \in \R^{n_2}_+ \times \R^{n_1}_+: \tilde{T}(\om) z + \tilde{W}(\om) y \leq \tilde{r}(\om)\}$. 
    Then, we have
    \begin{equation*}
    \begin{aligned}
    L(x, \lambda(\om), \om) 
        &
        = \max_{y \in \R^{n_2}_+, z \in \R^{n_1}_+} \Big\{q(\om)^\top y + x^\top G(\om) y + \lambda(\om)^\top z : \tilde{T}(\om) z + \tilde{W}(\om) y \leq \tilde{r}(\om) \Big\} \\
        &
        =\min_{\pi \in \R^{\tilde{m}}_+} \Big\{\tilde{r}(\om)^\top \pi : \tilde{T}(\om)^\top \pi \geq \lambda(\om), \tilde{W}(\om)^\top \pi \geq q(\om) + G(\om)^\top x
        \Big\},
    \end{aligned}
    \end{equation*}
    where $\pi$ represents dual variables.
    By incorporating the above into the Lagrangian dual~\cref{eq:dual_Q_om}, we have
    \begin{equation}
    \begin{aligned}
        D(x, \om)
        = \min_{\lambda(\om) \in \R^{n_1}, \pi \in \R^{\tilde{m}}_+} \Big\{\tilde{r}(\om)^\top \pi - \lambda(\om)^\top x : 
        \tilde{T}(\om)^\top \pi - \lambda(\om) \geq 0, \tilde{W}(\om)^\top \pi \geq q(\om) + G(\om)^\top x\Big\}.
    \end{aligned}
    \end{equation}
    Let $z \in \R_+^{n_1}$ and $y\in\R_+^{n_2}$ denote dual variables for the constraints. Then, the resulting dual formulation is
    \begin{equation}\label{eq:min-min_conv_form}
        \max_{y \in \R_+^{n_2}, z \in \R_+^{n_1}} \Big\{
            q(\om)^\top y + x^\top G(\om) y : \tilde{T}(\om) z + \tilde{W}(\om) y \leq \tilde{r}(\om), \ 
            z = x
        \Big\}.
    \end{equation}
    The feasible region of \cref{eq:min-min_conv_form} can be seen as the intersection of $\conv(Z(\om))$ and $\{(y, z): z = x\}$. By \Cref{assumption:binary_x_component_minmax}, $Z(\om)$ is independent of $z_i$ for $i \in [n_1] \setminus I$; thus, this intersection can be equivalently written as $\conv(Z(\om)) \cap \{(y, z): z_i = x_i, \forall i \in I\}$.
    Since $z_i$ is a binary variable for $i \in I$, the intersection defines a face of $\conv(Z(\om))$ for any $x \in X$. By the properties of faces, the $y$ component of any extreme point in this set satisfies the integrality constraints, i.e., $y \in \Y$. Consequently, the $y$ component of any optimal solution to \cref{eq:min-min_conv_form} is feasible to the recourse problem~\cref{eq:Q_om_re}, which implies that $D(x, \om) \leq Q(x, \om)$.
\end{proof}

\subsection{Proof of \texorpdfstring{\Cref{lem:sigma_optimal_condition}}{Lemma~\ref{lem:sigma_optimal_condition}}}
\begin{proof}{Proof}
    Fix $(\hat{x}, \hat{\lambda}(\om), \om) \in X \times \R_+^{n_1} \times \Om$ and let $(\hat{y}(\om), \hat{z}(\om))$ be an optimal solution of the Lagrangian relaxation~\cref{eq:lag_Q_om}. There can be two cases: \textit{(a)} $\hat{z}_i(\om) = \hat{x}_i$ for all $i \in I$ or \textit{(b)} $\hat{z}_i(\om) \neq \hat{x}_i$ for some $i \in I$.
    Clearly, if Case \textit{(a)}, then $\hat{y}(\om)$ is feasible to the recourse problem~\cref{eq:Q}, and thus $q(\om)^\top \hat{y}(\om) + \hat{x}^\top G(\om) \hat{y}(\om) \leq Q(\hat{x}, \om)$. Here, the equality holds since we have $Q(\hat{x}, \om) \leq L(\hat{x}, \hat{\lambda}(\om), \om) - \hat{\lambda}(\om)^\top \hat{x} = q(\om)^\top \hat{y}(\om) + \hat{x}^\top G(\om) \hat{y}(\om)$. This implies that the solution $\hat{y}(\om)$ is optimal to the recourse problem, and $\hat{\lambda}(\om)$ is an optimal solution to the Lagrangian dual.

    Now, consider Case \textit{(b)}. 
    Pick any $\bar{y} \in Y(\hat{x}, \om)$. Then, since solution $(y, z) = (\bar{y}, \hat{x})$ is non-optimal and feasible to the Lagrangian relaxation~\cref{eq:lag_Q_om}, we have
    \begin{equation}
        q(\om)^\top \hat{y}(\om) + \hat{x}^\top G(\om) \hat{y}(\om)  + \hat{\lambda}(\om)^\top \hat{z}(\om) > q(\om)^\top \bar{y} + \hat{x}^\top G(\om) \bar{y}  + \hat{\lambda}(\om)^\top \hat{x}.
    \end{equation}
    By rearranging the terms, we obtain
    \begin{equation}
        (q(\om) + G(\om)^\top \hat{x})^\top (\hat{y}(\om) - \bar{y}) > \hat{\lambda}(\om)^\top (\hat{x} - \hat{z}(\om)) 
        = \sum_{i \in I: \hat{z}_i(\om) \neq \hat{x}_i} \sigma_i
    \end{equation}
    Here, let $I_{d}$ denote the set $\{i \in I: \hat{z}_{i}(\om) \neq \hat{x}_i\}$. By \cref{assumption:complete_recourse}, $\max_{y \in Y(\hat{x}, \om)} (q(\om) + G(\om)^\top \hat{x})^\top (\hat{y}(\om) - y) < \infty$, and thus the left-hand side of the above inequality is bounded. This implies that there exists $\sigma_i \R_+, i \in I,$ such that the inequality is violated for all $\bar{y} \in Y(\hat{x}, \om)$.
    In other words, the following provides a sufficient condition under which it is enforced that $\hat{z}_i(\om) = \hat{x}_i$ for all $i \in I$:
    \begin{equation}
        \sum_{i \in I_d} \sigma_i
        \geq
        \max_{y \in Y(\hat{x}, \om)} 
        \Big\{
        (q(\om) + G(\om)^\top \hat{x})^\top (\hat{y}(\om) - y) 
        \Big\}.
    \end{equation}
\end{proof}

\subsection{Proof of \texorpdfstring{\Cref{thm:min-max_reform}}{Theorem~\ref{thm:min-max_reform}}}
\begin{proof}{Proof}
    By \Cref{lem:sigma_optimal_condition}, 
    in the reformulation~\cref{eq:min-max_bilinear}, we can fix $\lambda(\om)$ to its optimal value using the analytical form: $\lambda_{i}(\om) = \sigma_i (2 x_i - 1)$ for $i \in I$ and $\lambda_{i}(\om) = 0$ for $i \in [n_1]\setminus I$ for every $\om \in \Om$. Consequently, we can substitute the bilinear terms $\lambda(\om)^\top x$ and $\lambda(\om)^\top z$ in the objective functions as follows:
    $
        \lambda(\om)^\top x = \sum_{i \in I} \sigma_i (2 x_i - 1)x_i = \sum_{i \in I} \sigma_i x_i,
    $ 
    and 
    $
        \lambda(\om)^\top z = \sum_{i \in I} \sigma_i (2 x_i - 1) z_{i}
    $
    for $\om \in \Om$.
    This results in the formulation~\cref{eq:min-max_reform}.
\end{proof}

\subsection{Proof of \texorpdfstring{\Cref{prop:finite_convergence}}{Proposition~\ref{prop:finite_convergence}}}
\begin{proof}{Proof}
    Let $\phi(x)$ be the objective function of the reformulation~\cref{eq:min-max_reform}.
    To prove the statement, it suffices to show that \Cref{alg:min-max-L2} terminates with $UB=LB$ in a finite number of iterations, as this implies that $\phi(x^\ast) = UB = LB \leq \min_{x \in X} \phi(x) = \min_{x \in X} \{ f(x) + \sum_{\om \in \Om} p(\om) Q(x, \om) \}$. The last equality holds by \Cref{thm:min-max_reform}.
    
    For a solution $(y^l(\om), z^l(\om))$ and $\om \in \Om$, we have $\hat{Q}(x, \om) \geq q(\om)^\top y^l(\om) + x^\top G(\om) y^l(\om) + \sum_{i \in [n_1]}\sigma_i (2x_i - 1)z^l_i(\om)$ since the solution is feasible to the problem~\cref{eq:min-max_reform_Q}. Furthermore, the inequality is tight at $x^l$, i.e., $\hat{Q}(x^l, \om) = q(\om)^\top y^l(\om) + (x^l)^\top G(\om) y^l(\om) + \sum_{i \in [n_1]}\sigma_i (2x_i^l - 1)z^l_i(\om)$. Thus, the optimality cut generated using solutions $(y^l(\om), z^l(\om)), \om \in \Om,$ for iteration $l$ is a valid cut that supports the epigraph of $\sum_{\om \in \Om} p(\om) \hat{Q}(x, \om)$ at $x^l$, thus ensuring that the master problem exactly evaluates the objective value in the incumbent solution $x^l$. Since $|X| < \infty$, there exists $l < \infty$ such that $\hat{\theta} = \sum_{\om \in \Om} p(\om) \hat{Q}(x^{l+1}, \om)$ at \Cref{line:min-max_compute_lower_bound}. This implies that $UB = LB$, as $x^{l+1}$ is feasible to the original reformulation~\cref{eq:min-max_reform}.
 \end{proof}

\subsection{Proof of \texorpdfstring{\Cref{prop:min-min_valid_cut}}{Proposition~\ref{prop:min-min_valid_cut}}}
\begin{proof}{Proof}
    For scenario $\om$ and iteration $l$, the dual subproblem can be viewed as a restriction of the reformulated recourse problem~\cref{eq:min-min_reform_Q} where variables $\pi_{k}$ for $k=m^l + 1, \dots, m(\om)$ are restricted to be zero.
    Let $(\hat{\pi}(\om), z^l(\om))$ be a lifted solution of $(\pi^l(\om), z^l(\om))$, where $\hat{\pi}_{k}(\om) = \pi^l_{k}(\om)$ for $k=1, \dots, m(\om)$ and $\hat{\pi}_{k}(\om) = 0$ for $k=m^l + 1, \dots, m(\om)$. This solution is feasible to problem~\cref{eq:min-min_reform_Q}, and thus the following inequality holds:
    \begin{equation*}
    \begin{aligned}
        \hat{Q}(x, \om)
        & \geq \hat{\pi}(\om)^\top (\tilde{r}(\om) - \tilde{T}(\om) x) + \sum_{i \in [n_1]} \sigma_i z^l_{i}(\om) (2x_i - 1) \\
        & = \pi^l(\om)^\top (r^l(\om) - T^l(\om) x) + \sum_{i \in [n_1]} \sigma_i z^l_{i}(\om) (2x_i - 1).
    \end{aligned}
    \end{equation*}
    Multiplying the inequality by $p(\om)$ and aggregating them for all $\om \in \Om$ yields the optimality cut.
\end{proof}

\end{appendices}

\setstretch{1.25}
\bibliographystyle{apalike}
\bibliography{reference.bib}

\end{document}